\input amstex
\magnification=1095
\input amsppt.sty

\hsize=14truecm
\vsize=22.6truecm

\TagsOnLeft
\def\bino{\atopwithdelims[]}
\def\ch{\bold {char}}
\def\pa{\partial}
\def\lra{\longrightarrow}
\def\al{\alpha}
\def\be{\beta}
\def\ga{\gamma}
\def\si{\sigma}
\def\Th{\Theta}
\def\De{\Delta}
\def\de{\delta}
\def\th{\theta}
\def\vn{\varepsilon}
\def\vp{\varphi}
\def\ep{\epsilon}
\def\ot{\otimes}

\def\la{\lambda}
\def\La{\Lambda}

\rightheadtext{$q$-Divided Power Algebra, $q$-Derivatives and Quantum
Groups}
\leftheadtext{Naihong Hu}

\topmatter

\

\

\

\

\title{\bf Quantum Divided Power Algebra, $q$-Derivatives \\
and Some New Quantum Groups}
\endtitle
\author Naihong Hu \endauthor

\footnote " "{Supported in part by the
National Natural Science Foundation of China (Grant No. 19731004),
the Science Foundation of the University Doctoral Program CNCE,
the Ministry of Education of China and
the Shanghai Scientific and Technical Commission, as well as the Institute
of Mathematic of ECNU.}

\

\abstract
The discussions in the present paper arise from exploring intrinsically
the structure nature of the quantum $n$-space.
A kind of braided category $\Cal {GB}$ of $\La$-graded
$\th$-commutative associative algebras over a field $k$ is established.
The quantum divided power
algebra over $k$ related to the quantum $n$-space is introduced and
described as a braided Hopf algebra in $\Cal {GB}$
(in terms of its $2$-cocycle structure),
over which the so called special
$q$-derivatives are defined so that several new interesting quantum
groups, especially, the quantized polynomial algebra in $n$ variables
(as the quantized universal enveloping algebra of the abelian Lie algebra
of dimension $n$),
and the quantum group associated to the quantum $n$-space,
are derived from our approach independently of using the $R$-matrix.
As a verification of its validity
of our discussion, the quantum divided power algebra is equipped with
a structure of $U_q(\frak {sl}_n)$-module algebra via a certain $q$-differential
operators realization. Particularly, one of the four kinds of roots vectors
of $U_q(\frak {sl}_n)$ in the sense of Lusztig can be specified precisely
under the realization.

\endabstract

\endtopmatter
\document

\baselineskip13pt

\

\heading
{$\bold 1$. \ Introduction and Preliminaries}\endheading

\vskip0.12cm
\noindent
{\bf 1.1} \ It is a known fact that associated to a simple Lie group $G$
there exist two new Hopf
algebra structures, namely, quantum group $k_q[G]$, and quantized universal
enveloping algebra $U_q(\frak g)$. In particular, to $G=GL(n)$ or $SL(n)$ there
correspond two (right) comodule-algebras (cf. \cite{13}):
the first one
$$
k[A_q^{n|0}]=k\{x_1,\cdots,x_n\}/(x_jx_i-q\,x_ix_j, \ \; i<j)
$$
is the {\it quantum $n$-space} (i.e. {\it quantum symmetric algebra}) whereas
the second  one
$$
k[A_q^{0|n}]=k\{x_1\cdots,x_n\}/(x_i^2, \ x_j\,x_i+q^{-1}\,x_i\,x_j, \ \; i<j)
$$
is the {\it quantum exterior algebra} $\Lambda_q(n)$.

Similar to the classical case, there exists a Hopf duality between
$k_q[G]$ and  $U_q(\frak g)$ (cf. \cite{2}, \cite{6}, \cite{7}, etc.). Hence,
$k[A_q^{n|0}]$ and $k[A_q^{0|n}]$ should be (left) $U_q(\frak g)$-module
algebras (cf. \cite{1, 16}), where $\frak g=\frak {gl}_n$, $\frak {sl}_n$.  Thus,
a natural problem arises from here, i.e. how to concretely realize the quantized
universal enveloping algebra $U_q(\frak g)$ as certain $q$-differential
operators over the associated quantum $n$-space such that $k[A_q^{n|0}]$
becomes a $U_q({\frak g})$-module algebra. According to \cite{12}, this
question has its (non-commutative) geometric meaning.

\vskip0.12cm
\noindent
{\bf 1.2} \ For the case when $\frak g$ is $\frak {sl}_2$, this was solved
independently by J. Wess and B. Zumino (\cite{17}), S. Montgomery and S.P.
Smith (in fact, their realization is
just relative to the Woronowicz's Hopf subalgebra of $U_q(\frak {sl}_2)$)
(\cite{12}), and C. Kassel (\cite{7}). Kassel's treatment depends on
complicated commutative operations concerning
left- and right- multiplications so that it seems impossible to
generalize to the general case $\bold A_{n-1}$ by following his approach
(Note that in his realization, the two $q$-derivatives
$\partial_q/\partial x$, $\partial_q/\partial y$ are commutative).
For $\frak {sl}_n$,
T. Hayashi (\cite{5}) gave a realization of $U_q(\frak {sl}_n)$ over
a polynomial ring $S=k[X_1,\cdots,X_{n+1}]$ (also see Jantzen's
book \cite{6}), however, his realization cannot make the polynomial ring
$S$ a $U_q(\frak {sl}_n)$-module algebra in the sense of \cite{1, 16}.

In the present paper, we first establish a kind of braided
category $\Cal {GB}$ of $\La$-graded $\th$-commutative algebras over a field
$k$ in section 2, where $\La$ is a free abelian group (of finite rank) and
$\th$ is a bicharacter (or called a $2$-cocycle) on $\La$, and then
describe $k[A_q^{n|0}]$ as a
braided Hopf algebra in the category $\Cal {GB}$,
which is a generalization of super-Hopf algebra
(cf. \cite{15}) or a kind of braided Hopf algebra in the sense of S. Majid
(cf. \cite{10}), relative to a 2-cocycle defined on the $\Bbb Z^n$-graded
structure of the quantum $n$-space.
Instead of it, we introduce the notion of
quantum divided power algebra we would work with.
In section 3, we define a kind of
$q$-differential operators (i.e. special $q$-derivatives)
over it which differs from the above-mentioned and
can be used to give the required realization of $U_q(\frak g)$ as in {\bf 1.1}
(for $\frak g=\frak {gl}_n$, or $\frak {sl}_n$). Also, we obtain several
new quantum groups, for instance, the
quantum group $\frak D_q$ whose smash product relative to the quantum
divided power algebra $\Cal A_q$ yields a
kind of quantum Weyl algebra (i.e. the algebra of quantum differential
operators) distinguished from those appeared in the literature (for instance,
\cite{3}, \cite{5}, \cite{10}, etc.) to the best of my knowledge. Of
particular interest in our discussion, we are able to obtain the exact
object of polynomial algebra in $n$ variables in the context of quantum
groups as well as the quantum group structure associated to the quantum
$n$-space $k[A_q^{n|0}]$.
On the other hand, with the realization in section 4, we consider the submodules
structure of the quantum divided power (restricted) algebra, especially,
we show that one of the four kinds of roots vectors of $U_q(\frak {sl}_n)$
introduced by G. Lusztig in \cite{9} can be specified precisely with those
$q$-differential operators we defined.

\vskip0.12cm
For the sake of latter discussion, we recall the following notions.

\vskip0.12cm
\noindent
{\bf 1.3} \ Recall that a Hopf algebra $(H,\Delta,\ep,S)$ over $k$
means $H$ is an algebra,
$\Delta: H\to H\otimes H$ (the comultiplication) and
$\ep: H\to k$ (the counit) are algebra homomorphisms and $S: H\to H$ (the
antipode) plays the role of inverse. Here $H\otimes H$ is of the tensor
product algebra structure. Call an
(associative) algebra $A$ over $k$ an {\it $H$-module algebra} (cf. \cite{1},
\cite{16}, etc.) if $A$ has an
(left) $H$-module structure such that
$$
\gather
h\,1_A =\epsilon(h)\,1_A,\tag{i}\\
h(a\,b) =\sum (h_{(1)}\,a)\,(h_{(2)}\,b),\tag{ii}
\endgather
$$
for $h\in H$, $a, b\in A$ with $\Delta(h)=\sum h_{(1)}\ot h_{(2)}$. Here the
second condition means that the multiplication is a homomorphism of
$H$-modules.

Given two automorphisms $\sigma$ and $\tau$ of an algebra $A$, a linear
endomorphism $\delta$ of $A$ is called a {\it $(\sigma,\tau)$-derivation} if
$$
\delta(a\,a')=\sigma(a)\,\delta(a')+\delta(a)\,\tau(a'),\leqno(\text{\rm iii})
$$
for $a, a'\in A$ (cf. \cite{7}).

A {\it quantum group} in the sense of Drinfeld (cf. \cite{4})
is a non-commutative and non-cocommutative Hopf algebra.
For super-quantum groups and super-Hopf algebras, however,
the main difference is that the algebra structure on $H\otimes H$ uses
the super-transposition (cf. \cite{15}): $\psi(x\otimes y)
=(\pm 1)^{|x||y|}y\otimes x$ (on homogeneous elements). More generally,
by a {\it braided Hopf algebra} in a certain {\it braided category} means the algebra
structure on $H\otimes H$ is provided by a certain ``braiding'' in the
sense of Majid (cf. \cite{10}). Precisely, between any two objects
there is a tensor product that is commutative
and associative up to isomorphism. The first of these isomorphisms is the
braiding $\psi_{V,W}:V\otimes W\to W\otimes V$. For any objects $U, V$ and
$W$, the braiding obeys
$$
\gather
\psi_{U\otimes  V,W}=(\psi_{U,W}\ot\text{id}_V)\circ
(\text{id}_U\ot\psi_{V,W}),\tag{\text{\rm iv}} \\
\psi_{U,V\otimes W}
=(\text{id}_V\ot\psi_{U,W})\circ(\psi_{U,V}\ot\text{id}_W).
\endgather
$$
In addition, there should be an identity object $\underline 1$ for $\otimes$,
and one has $\psi_{\underline 1,V}=\text{id}=\psi_{V,\underline 1}$.
In these formulae, the associativity isomorphism and isomorphisms such as
$V\otimes\underline 1\cong V\cong \underline 1\otimes V$ are suppressed.

\vskip0.12cm
\noindent
{\bf 1.4} \ The {\it quantized universal enveloping algebra}
$U_q(\frak {sl}_n)$ (cf. \cite{6}, \cite{7}, \cite{9}, etc.) is the $k$-algebra
generated by the symbols $\Cal K_i^{\pm 1}$,
$e_i$ and $f_i$ $(1\le i\le n-1)$ with the following defining relations:
$$
\gather
 \Cal K_i\,\Cal K_i^{-1}=\Cal K_i^{-1}\,\Cal K_i=1, \qquad
 \Cal K_i\,\Cal K_j=\Cal K_j\,\Cal K_i,\tag{i}\\
 \Cal K_i\,e_j\,\Cal K_i^{-1}=q^{a_{ij}}e_j,
 \qquad \Cal K_i\,f_j\,\Cal K_i^{-1}=q^{-a_{ij}}f_j,\tag{ii}\\
 [e_i,f_j]=\delta_{ij}\frac{\Cal K_i-\Cal K_i^{-1}}{q-q^{-1}},\tag{iii}\\
 e_i^2\,e_j-(q+q^{-1})\,e_i\,e_j\,e_i+e_j\,e_i^2=0 \qquad (|i-j|=1),\tag{iv}\\
 e_i\,e_j=e_j\,e_i \qquad (|i-j|>1),\\
 f_i^2\,f_j-(q+q^{-1})\,f_i\,f_j\,f_i+f_j\,f_i^2=0 \qquad (|i-j|=1),\tag{v}\\
 f_i\,f_j=f_j\,f_i \qquad (|i-j|>1),
\endgather
$$
where $q\in k^*$ and $(a_{ij})$ is the Cartan matrix of type $\bold A_{n-1}$.

The Hopf algebra structure of $U_q(\frak {sl}_n)$ is as follows:
$$
\gather
 \Delta(\Cal K_i^{\pm 1})=\Cal K_i^{\pm 1}\ot\Cal K_i^{\pm 1},\tag{vi}\\
 \epsilon(\Cal K_i^{\pm 1})=1,\\
 S(\Cal K_i^{\pm 1})=\Cal K_i^{\mp 1},\\
 \Delta(e_i)=e_i\ot\Cal K_i+1\ot e_i,\tag{vii}\\
 \epsilon(e_i)=0,\\
 S(e_i)=-e_i\,\Cal K_i^{-1},\\
 \Delta(f_i)=f_i\ot 1+\Cal K_i^{-1}\ot f_i,\tag{viii}\\
 \epsilon(f_i)=0,\\
 S(f_i)=-\Cal K_i\,f_i.
\endgather
$$

Let $P$ be the weight lattice for $\frak {gl}_n$. It is the free
$\Bbb Z$-module of rank $n$ with canonical basis
$\{\vn_i\}_{1\le i\le n}$, and
$\{\alpha_i=\vn_i-\vn_{i+1}\mid 1\le i\le n-1\}$
is the set of simple roots of $\frak {gl}_n$.
$Q=\bigoplus_{i=1}^{n-1}\Bbb Z\alpha_i\subset P$ is the root lattice of
$\frak {gl}_n$. Fix a symmetric bilinear form $\langle ,\rangle:
P\times P\lra\Bbb Z$ such that $\langle\vn_i,\vn_j\rangle=
\delta_{ij}$ for $1\le i,j\le n$. Through this pairing, we can identify $P$
with its dual $P^*=\text{Hom}_{\Bbb Z}(P,\Bbb Z)$.

Now we can state the presentation of $U_q(\frak {gl}_n)$ as follows.
Change (i), (ii) into
$$
\gather
 k_i\,k_i^{-1}=k_i^{-1}\,k_i=1, \qquad
 k_i\,k_j=k_j\,k_i \qquad (1\le i, j\le n),\tag{ix}\\
 \Cal K_i=k_i\,k_{i+1}^{-1} \qquad (1\le i\le n-1),\\
 k_i\,e_j\,k_i^{-1}=q^{\langle\vn_i,\alpha_j\rangle}e_j, \qquad
 k_i\,f_j\,k_i^{-1}=q^{-\langle\vn_i,\alpha_j\rangle}f_j,\tag{x}
\endgather
$$
but keep (iii)--(v) invariant. As for its Hopf algebra structure,
we only need to replace $\Cal K_i$ $(1\le i\le n-1$) in (vi) with
$k_i$ $(1\le i\le n)$, except with the same items as (vii) and
(viii).

\medskip
\noindent
{\bf 1.5} \
As we know, the $q$-binomial coefficients are closely related to the
study of quantizations $U_q(\frak g)$ of enveloping algebras. For our work
here we need some known facts involving them.

Let $\Bbb Z[v,v^{-1}]$ be the Laurent polynmmial ring in variable
$v$. For any integer $n\ge 0$ we define $[n]_v=\frac{v^n-v^{-n}}{v-v^{-1}}
\in\Bbb Z[v,v^{-1}]$, and $[n]_v!=[n]_v[n-1]_v\cdots [1]_v$. It is well
known that for two integers $m$, $r$ with $r\ge 0$ we have (cf. \cite{8}):
$$
{m\bino r}_v
=\prod_{i=1}^r\frac{v^{m-i+1}-v^{-m+i-1}}{v^i-v^{-i}}\in\Bbb Z[v,v^{-1}].
\leqno(\text{\rm i})
$$
By (i), we get
${m\bino r}_v=\frac{[m]_v!}{[r]_v!\,[m-r]_v!}$ if $0\le r\le m$,
${m\bino r}_v=0$ if $0\le m<r$,
and ${m\bino r}_v=(-1)^r\,{-m+r-1\bino r}_v$ if $m<0$.
We again set ${m\bino r}_v=0$ when $r<0$.

Suppose now that $k$ is a field and $q\in k^*$.
By definition, we get $[n]:=[n]_{v=q}$,
$[n]!:=[n]_{v=q}!\in k$ for $n\ge 0$ and ${m\bino r}:={m\bino r}_{v=q}\in k$
when $v$ is specialized to be $q$.
Note that the $q$-binomial coefficients ${n\bino r}$ ($0\le r\le n$)
can be defined recursively by
$$
{n\bino r}=q^{r-n}\,{n-1\bino r-1}+q^r\,{n-1\bino r}, \qquad
{1\bino 0}=1={1\bino 1}.$$

The combinatorial formula below involving the $q$-binomial
coefficients ${n\bino r}$ is well known (cf. \cite{8}).
$$
\prod_{i=0}^{n-1}(1+q^{2i}\,x)=\sum_{r=0}^nq^{(n-1)r}{n\bino r}\,x^r.
\leqno(\text{\rm ii})
$$

The situation when $q$ is a primitive root of unity is of particular
interest in quantum phenomenon.
We here introduce a notion, so called the {\it characteristic} of $q$,
which is defined as the minimal positive integer $l$ such that $[\,l\,]=0$
and denoted by $\ch(q)=l$. Also we define $\ch(q)=0$ when $q$ is generic
(in this case, $[\,n\,]\ne 0$ for any non-zero $n\in\Bbb Z$).
If $q\ne\pm 1$, then $\ch(q)=l\,(>2)$ implies two cases
that $q$ is a primitive $2l$-th root of $1$, or a primitive $l$-th root
of $1$ but $l$ odd. And vice versa.

The following Lemma is proved by Lusztig under the assumption
$l\ge 3$ is odd.
\proclaim{Lemma} \ Assume that $q\in k^*$ and
$\ch(q)=l\ge 3$.

$(1)$
Let $m=m_0+m_1l$, $r=r_0+r_1l$ with $0\le m_0$, $r_0<l$, $m_1$, $r_1\ge 0$
and $m\ge r$. Then ${m\brack r}=
{m_0\brack r_0}\binom{m_1}{r_1}$, where $\binom{m_1}{r_1}$ is an ordinary
binomial coefficient.

$(2)$ Let $m=m_0+m_1l$, $0\le m_0<l$, $m_1\in\Bbb Z$. Then
${m\bino l}=m_1$ if $m_1\ge 0$; and ${m\bino l}=-(-1)^lm_1$ if $m_1<0$ $($i.e.
$m<0)$.

$(3)$ If $m=m_0+m_1l$, $m'=m_0'+m_1'l\in\Bbb Z$ with $0\le m_0$, $m_0'<l$
satisfy
$q^m=q^{m'}$, ${m\bino l}={m'\bino l}$,
then $m=m'$ if $l$ is odd or $l$ even but $m_1m_1'\ge 0$; and
$m'=\overline m=m_0-m_1l$ if $l$ is even but $m_1m_1'<0$.
\endproclaim
\demo{Proof} By using formula (ii) and following the same argument as
the proofs of Proposition 3.2 \& Corollary 3.3
in \cite{8}, we readily show that (1),  2) and the first claim of (3) hold.
For the case where $l$ is even (in this case, $q$ must be a primitive
$2l$-th root of $1$) and $m_1m_1'<0$, we deduce from
$q^m=q^{m'}$ that $m_0=m_0'$ and $m_1=m_1'+2r$. Since $m_1m_1'<0$, we
can let $m_1>0$ and $m_1'<0$. By (2), we
get $m_1={m\bino l}={m'\bino l}=-m_1'$, as required.
\hfill\qed
\enddemo

\

\vskip0.35cm
\heading{$\bold 2$. \ Braided Hopf Algebra and Quantum Divided Power Algebra}
\endheading
\vskip0.12cm

\noindent
{\bf 2.1} \ Let $\alpha=(\al_1,\cdots,\al_n), \be=(\be_1,\cdots,\be_n)
\in \Bbb Z^n$ be any two integers $n$-tuples, and define a product of them by
$$
 \al* \be=\sum_{j=1}^{n-1}\sum_{i>j}\al_i\be_j,\leqno (\text{\rm i})
$$
then from (i), one has
\proclaim
{Lemma} \
The product $*$ satisfies the following distributive laws:

$(1)$ \ $( \al+ \be)* \ga= \al * \ga+ \be * \ga$,

\qquad $\al *( \be+ \ga)= \al * \be+ \al * \ga$, \quad
in particular,

$(2)$ \ $\vn_i*\be
=\sum_{s<i}\be_s$
\quad $(1\le i\le n)$,

\qquad $(\vn_i-\vn_{i+1})* \be
=-\be_i$
\quad $(1\le i <n)$.

$(3)$ \ $\be*\vn_i=\sum_{s>i}\be_s$
\quad $(1\le i\le n)$,

\qquad$ \be *(\vn_i-\vn_{i+1})=\be_{i+1}$
\quad $(1\le i <n)$.

\noindent
Here $\vn_i=(\delta_{1i},\cdots,\delta_{ni})$ $(1\le i\le n)$ is a
basis of $\Bbb Z^n$ as $\Bbb Z$-module.
\hfill\qed
\endproclaim

Now for $\al, \be\in \Bbb Z^n$ and $q\in k^*$, we define a mapping
$\theta: \Bbb Z^n\times \Bbb Z^n\to k^*$ by
$$
\theta(\al,\be)=q^{\al *\be-\be *\al}.\leqno(\text{\rm ii})
$$
In particular,
$$\th(\vn_i,\vn_j)=\cases
q, & i>j,\\
1, & i=j,\\
q^{-1}, & i<j.
\endcases\leqno(\text{\rm iii})
$$

Obviously, the mapping $\theta$
has the following properties:
$$\gather
\th(\al+\be,\ga)=\th(\al,\ga)\th(\be,\ga),\tag{iv}\\
\th(\al,\be+\ga)=\th(\al,\be)\th(\al,\ga),\tag{v}\\
\th(\al,0)=1=\th(0,\al),\tag{vi}\\
\th(\al,\be)\th(\be,\al)=1=\th(\al,\al).\tag{vii}
\endgather
$$
Actually, such a mapping $\th$ with the above properties is
a {\it bicharacter} of the additive group $\Bbb Z^n$.

Let $\al=(\al_1,\cdots,\al_n)\in\Bbb Z_+^n$ be any
nonnegative-integers $n$-tuple, $x^{\al}=x_1^{\al_1}\cdots
x_n^{\al_n}$ be any nonzero monomial in $k[A_q^{n|0}]$, then
$\{x^{\al}\mid \al\in\Bbb Z_+^n\}$ constitutes a canonical basis of
$k[A_q^{n|0}]$. Thus by definition (see {\bf 1.1}),
$k[A_q^{n|0}]=\bigoplus_{\al\in\Bbb Z^n} kx^\al$ is a $\Bbb
Z^n$-graded algebra (with $x^\al=0$ for $\al\not\in\Bbb Z_+^n$). Set
$x^{\vn_i}=x_i$.

\medskip
\noindent
{\it Remark.} \ We will point out that the quantum $n$-space $k[A_q^{n|0}]$
has a so-called ``braided'' Hopf algebraic structure (with respect to the
above defined bicharacter $\th$). Actually, we can show a more general
fact in the next subsection.

\medskip
\noindent
{\bf 2.2} \ Let $\Lambda$ be a (finitely generated) abelian group and
$\th$ a bicharacter of the abelian group $\Lambda$ (namely, $\th$ is of
properties {\bf 2.1} (iv) --- (vii)).
Recall a {\it $2$-cocycle} $\eta$ defined on $\Lambda$ with coefficients
in $k^*$, {\it i.e.,} a mapping
$\eta: \Lambda\times\Lambda\to k^*$
satisfying the relation
$$
\eta(\al,\be)\,\eta(\al+\be,\ga)=\eta(\be,\ga)\,\eta(\al,\be+\ga).
\leqno(\text{\rm i})
$$

Observe that any bicharacter $\th$ of $\Lambda$ is a $2$-cocycle
on $\Lambda$, which is useful to our next consideration.

Let $\Cal {GA}$ denote a category of graded associative unitary
algebras over $k$. That is, for any $(A, \La_A)\in\text{\it Ob}\,(\Cal {GA})$,
$A=\bigoplus_{\al\in\La_A}A_\al$ is a $\La_A$-graded associative algebra
with $k\subseteq A_0$ and $A_\al\cdot A_\be\subseteq A_{\al+\be}$, where
$\La_A$ is an abelian group.
Moreover, $\phi_{A,B}:=(\phi,\vp): (A, \La_A)\to (B, \La_B)$ is a morphism
between $(A,\La_A)$ and $(B,\La_B)$, if
$\phi: A\to B$ is a graded algebra homomorphism
and $\vp: \La_A\to\La_B$ is a group homomorphism such that
$\phi(A_\al)\subseteq B_{\vp(\al)}$.

Now for any object $(A,\La)\in\text{\it Ob}\,(\Cal {GA})$, to
$\th$ an arbitrary bicharacter of $\La$, we can associate an
opposite object $(A^{\text{op}},\La)\in\text{\it Ob}\,(\Cal {GA})$ as follows:

Denote $A^{\text{op}}:=(A,\circ)$, where
$A^{\text{op}}=\bigoplus_{\al\in\La}A_\al$.
Define
$$
a\circ b:=\th(\al,\be)\,ba, \qquad \forall \ a\in A_\al, \ b\in A_\be.
\leqno(\text{ii})
$$
Clearly, $A_\al\circ A_\be\subseteq A_{\al+\be}$, for any
$\al$, $\be\in\La$. On the other hand, by the $2$-cocycle
property (i) of $\th$ as a bicharacter, for any $a\in A_\al$, $b\in A_\be$
and $c\in A_\ga$, we have
$$
\split
(a\circ b)\circ c&=\th(\al,\be)\,\th(\al+\be,\ga)\,c(ba)\\
&=\th(\be,\ga)\,\th(\al,\be+\ga)\,(cb)a=a\circ (b\circ c).
\endsplit
$$
This means $A^{\text{op}}$ is an associative algebra. Thus
$(A^{\text{op}},\La)\in\text{\it Ob}\,(\Cal {GA})$.

Consider a subcategory $\Cal {GC}$ in category $\Cal {GA}$,
whose objects consist of $\Lambda$-graded $\th$-commutative algebras over $k$
(where $\La$ is an arbitrary abelian group, $\th$ is an arbitrary
bicharacter of $\La$),
{\it i.e.}, any object $(A, \La_A; \th_A)$ in $\Cal {GC}$
is called {\it $\La_A$-graded $\th_A$-commutative}, if
$(A,\La_A)\in\text{\it Ob}\,(\Cal {GA})$ and
$x\cdot y=\th_A(\al,\be)\,y\cdot x$, $\forall\; x\in A_\al$, $y\in A_\be$.
Moreover, $\Phi_{A,B}:=(\phi_{A,B},\tilde\vp): (A,\La_A;\th_A)\to
(B,\La_B;\th_B)$ is a morphism,
if $\phi_{A,B}: (A,\La_A)\to (B,\La_B)$
is a morphism in category $\Cal {GA}$ such that
$\tilde \vp(\th_A)=\th_B$, {\it i.e.},
$\th_A(\al,\be)=\th_B(\vp(\al),\vp(\be))$.

For any $(A,\La,\th)\in\text{\it Ob}\,(\Cal{GC})$, by (ii), we notice that
$a\circ b=\th(\al,\be)\,ba=ab$, for any $a\in A_\al$, $b\in A_\be$.
This means the opposite of any object in the category $\Cal {GC}$
coincides with the itself, which is the case we have in the category
of (usual) commutative algebras.

\medskip
\noindent
{\it Examples.} \ Any commutative algebra over $k$ in the usual sense can be
considered an object in $\Cal {GC}$ with a trivial grading relative
to a trivial group. The polynomial algebra $k[t_1,\cdots,t_n]$
is a $\Bbb Z^n$-graded $\th$-commutative algebra with
a trivial bicharacter $\th$ (means $\th(\al,\be)\equiv 1$  for any
$\al$, $\be\in \Bbb Z^n$). Again, the quantum $n$-space
$k[A_q^{n|0}]=\bigoplus_{\al\in\Bbb Z^n} kx^\al$
is a $\Bbb Z^n$-graded $\th$-commutative algebra with the
bicharacter $\th$ of $\Bbb Z^n$ defined in {\bf 2.1} (ii).

\medskip
Now assume that $\La$ is a free abelian group and
$\th: \La\times \La\to k^*$ is a non-trivial bicharacter of $\La$.
Suppose that $F=\bigoplus_{\al\in\La}F_\al\in \text{\it Ob}\,(\Cal {GA})$
is a free object, that is, $F$ is a free $\La$-graded associative
algebra (with $1$).

Let $I=\langle\, x\cdot y-\th(\al,\be)\,y\cdot x,
\ \forall\; x\in F_\al, \ y\in F_\be, \ \forall\;\al, \ \be\in \La\,\rangle$
denote the $\La$-graded ideal generated by all homogeneous elements of form
$x\cdot y-\th(\al,\be)\,y\cdot x$, $x\in F_\al$, $y\in F_\be$.
Set $\Cal F:=F/I$. Thus $\Cal F$ is a free $\La$-graded $\th$-commutative
associative algebra over $k$, that is, $\Cal F\in \text{\it Ob}\,(\Cal {GC})$.

\medskip
\noindent {\it Remark.} \ In such $\La$-graded $\th$-commutative
associative algebras, ``$\th$-commutative'' is well-defined due to
the properties {\bf 2.1} (iv)---(vii) of the bicharacter $\th$, on
the other hand, the $2$-cocycle property (i) of the $\th$ ensures
the compatibility between ``$\th$-commutativity'' and
``associativity''.

\medskip\noindent
{\bf 2.3} \
Fix an abelian group $\La$ and a bicharacter
$\th$ of $\La$, we consider a subcategory $\Cal {GB}$
in category $\Cal {GC}$
relative to the pair $(\La,\th)$, where
$\forall\; (A,\La;\th)$, $(B,\La;\th)\in\text{\it Ob}\,(\Cal {GC})$,
the morphisms between $(A,\La;\th)$ and $(B,\La;\th)$ in
$\Cal {GC}$ are of forms:
$\Phi_{A,B}=(\phi_{A,B},\text{id})=((\phi,\text{id}),\text{id}\,)\equiv \phi:
(A,\La;\th)\to (B,\La;\th)$, where $\phi: A\to B$ is a graded
algebra homomorphism such that $\phi(A_\al)\subseteq B_\al$ for all
$\al\in\La$.

Note that $\Cal {GB}$ is a {\it braided category} in the sense of section
{\bf 1.3}.
First of all, there exists an identity object $\underline 1=k=
\bigoplus_{\al\in\La}k_\al\in\text{\it Ob}\,(\Cal {GB})$ with
$k_\al=\delta_{0,\al}\,k$.
Next, for any
$U=\bigoplus_{\al\in\Lambda}U_\al$,
$V=\bigoplus_{\be\in\Lambda}V_\be\in \text{\it Ob}\,(\Cal{GB})$, we have
$$
U\otimes V=\bigoplus_\ga\;(U\ot V)_\ga
=\bigoplus_{\ga=\al+\be}\,U_\al\otimes V_\be\in \text{\it Ob}\,(\Cal {GB}).
$$
Define the mapping $\psi_{U,V}: U\otimes V \to V\otimes U$ as
$$
\psi(x\otimes y)=\th(\al,\be)\; y\otimes x, \quad \text{\it for\; } \; x\in
U_\al, \; y\in V_\be.\leqno(\text{i})
$$
The properties of a bicharacter (see {\bf 2.1} (iv)---(vii)\,)
ensure that the mapping $\psi$ is a braiding in the
sense of {\bf 1.3} (iv).
In particular, for any algebra $H\in\Cal {GB}$, we have
its opposite object, $H^{\text{op}}\equiv H\in\text{\it Ob}\,(\Cal {GB})$, and
its tensor object, $H\ot H\in\text{\it Ob}\,(\Cal {GB})$,
whose algebra structure is given by
$$
(a\ot b)(c\ot d)=a\psi(b\ot c)d=\th(\al,\be)\;ac\ot bd, \quad
\text{\it for} \; a, d\in H, b\in H_\al, c\in H_\be.\leqno(\text{ii})
$$

Let $u(H)$ denote the group of invertible elements of $H$.

\proclaim
{Theorem} \  Let $\La$ be a free abelian group with
a basis $\{\al_1,\cdots,\al_n\}$,
$\th$ a bicharacter
of $\La$ and $\Cal {GB}$ the braided category relative to the pair $(\La, \th)$.
Suppose that $H=\bigoplus_{\al\in\Lambda}H_\al\in
\text{\it Ob}\,(\Cal {GB})$ such that $u(H)=k^*$.
If $H$ has no zero divisors $\ne 0$
and $\{\,a_{ij_i}\in H_{\al_i}\mid 1\le j_i\le s_i, 1\le i\le n\,\}$
is a set of generators for $H$.
Then $(H, \Delta, \ep, S)$
is a braided-commutative Hopf algebra relative to a braiding
$\psi: H\otimes H\to H\otimes H$ defined by $\psi(a\otimes b)
=\th(\al,\be)\; b\otimes a$ for $a\in H_\al$, $b\in H_\be$,
where the mappings $\Delta$, $\ep$ and $S$ defined below
$$\gather
\Delta: H\to H\otimes H, \qquad \Delta(a)=a\otimes 1+1\otimes a, \quad
\text{\it for} \; a\in H_{\al_i},\\
\ep: H\to k, \qquad \ep(a)=\delta_{0,\al_i}a,
\quad \text{\it for} \; a\in H_{\al_i},\\
S: H\to H^{\text{op}}, \qquad S(a)=-a, \quad \text{\it for} \; a\in H_{\al_i}
\endgather
$$
are the morphisms in the braided category $\Cal {GB}$ $($where
$S(ab)=\th(\al,\be)\,S(b)S(a)=S(a)\circ S(b)$, $\forall \ a\in H_\al$, $b\in H_\be$ $)$.

In particular, any free object $\Cal F$ in $\Cal {GB}$ constructed in
{\bf 2.2} is a braided-commutative Hopf algebra.
\endproclaim
\demo{Proof}
First we need to check that $\Delta$, $\ep$ and $S$ preserve
the algebraic relations of $H$.

For $a\in H_{\al_i}$, $b\in H_{\al_j}$, using (ii), we have
$$
\split
\Delta(a)\Delta(b)&=(a\ot 1+1\ot a)(b\ot 1+1\ot b)\\
&=ab\ot 1+\th(\al_i,\al_j)b\ot a+a\ot b+1\ot ab\\
&=\th(\al_i,\al_j)(ba\ot 1+b\ot a+\th(\al_j,\al_i)a\ot b+1\ot ba)\\
&=\th(\al_i,\al_j)\Delta(b)\Delta(a),\\
\ep(ab)&=\delta_{0,\al_i}\delta_{0,\al_j}ab=0=\th(\al_i,\al_j)\ep(ba),\\
S(ab)&=\th(\al_i,\al_j)\,ba=ab=S(\th(\al_i,\al_j)\,ba).
\endsplit
$$

Based on the consideration above, as well as the actions of $\Delta$ and
$\ep$ on the generators of $H$, it is readily to see that
$(1\ot\Delta)\Delta=(\Delta\ot 1)\Delta$ and
$(1\ot \ep)\Delta=1=(\ep\ot 1)\Delta$ hold.
Also, since $S(a)+S(1)a=0$, $\forall \, a\in H_{\al_i}$, we get
$m\circ (S\ot 1)\circ\Delta=\eta\circ\ep=m\circ (1\ot S)\circ\Delta$
(where $(H,m,\eta)$ is the algebra structure of $H$).
Hence, $(H,m,\eta,\Delta,\ep,S)$ is a braided-commutative Hopf algebra.
\hfill\qed
\enddemo

\proclaim
{Corollary} The quantum $n$-space $k[A_q^{n|0}]$ with respect to
the bicharacter $\th$ of $\Bbb Z^n$
given in {\bf 2.1} $(\text{\rm ii})$ is a $\Bbb Z^n$-graded
braided-commutative Hopf algebra with $\Delta(x_i)=x_i\ot 1+1\ot x_i$,
$S(x_i)=-x_i$, and $\ep(x_i)=0$ for $x_i\in k[A_q^{n|0}]_{\vn_i}$, and
where the braiding $\psi(a\ot b)=\th(\al,\be)\;b\ot a$ for $a\in k[A_q^{n|0}]_\al$,
$b\in k[A_q^{n|0}]_\be$.
\hfill\qed
\endproclaim

\noindent
{\bf 2.4} \ Here we will equip the quantum $n$-space $k[A_q^{n|0}]$
with a divided power structure when $\ch(q)=0$.
More generally, we can introduce
a quantum divided power algebra $\Cal A_q(n)$ for an arbitrary $q\in k^*$
as follows.

Let $\Cal A_q(n):=\langle x^{(\al)}\mid \al\in\Bbb Z_+^n\rangle$
be a vector space over $k$ generated by the (divided power) basis
$x^{(\al)}$ ($\al\in\Bbb Z_+^n$) with $x^{(0)}=1$.
Define the multiplication in $\Cal A_q(n)$ by
$$
x^{(\al)}\,x^{(\be)}=q^{\al*\be}\,{{\al+\be}\bino{\al}}\,x^{(\al+\be)}
=\th(\al,\be)\,x^{(\be)}\,x^{(\al)},
\leqno(\text{i})
$$
where
${{\al+\be}\bino{\al}}:=\prod_{i=1}^n{{\al_i+\be_i}\bino{\al_i}}$  and
${{\al_i+\be_i}\bino{\al_i}}=\frac{[\al_i+\be_i]!}{[\al_i]![\be_i]!}$
for $\al_i$, $\be_i\in\Bbb Z_+$.
By slight abuse of notation, we also write $x_i=x^{(\vn_i)}$.

Obviously, $\Cal A_q(n)$ with the above multiplication $\text{(i)}$
forms an associative algebra, and we call it a {\it quantum divided
power algebra}.

In particular, when $\ch(q)=0$, for $x^\al\in k[A_q^{n|0}]$,
we set $x^{(\al)}:=\frac{1}{[\al]!}\,x^\al$ for $\al\in\Bbb Z_+^n$,
where $[\,\al\,]!:=\prod_{i=1}^n[\,\al_i\,]!$.
According to the multiplication of $k[A_q^{n|0}]$, we see that
the algebraic structure $k[A_q^{n|0}]$ coincides with $\Cal A_q(n)$, and
$\{x^{(\al)}\mid \al\in\Bbb Z_+^n\}$ forms a (divided power)
basis of $k[A_q^{n|0}]$ with $x^{(0)}=1$, which gives the
{\it divided power structure} over the quantum $n$-space $k[A_q^{n|0}]$.
In the case, $\{x_i\mid 1\le i\le n\}$ is a set of generators of $\Cal A_q$.

Of particular interest to us in introducing $\Cal A_q(n)$ is the case
when $\ch(q)=l\ge 3$, in which we can obtain a finite dimensional
quantum divided power algebra.

Set $\tau=(l-1,\cdots,l-1)$. By $\al\le \be$ means $\al_i\le \be_i$
for all $i$. Denote $\Cal A_q(n,\bold 1):=\langle x^{(\al)}\mid
\al\in\Bbb Z_+^n, \al\le \tau\rangle$. Note that $[\,l\,]=0$.
Especially for $0\le s, t< l$, we have ${{s+t}\bino{s}}=0$ if
$s+t\ge l$. Consequently, the subspace $\Cal A_q(n,\bold 1)$ is
closed under the multiplication $\text{(i)}$ in $\Cal A_q(n)$,
namely, $\Cal A_q(n,\bold 1)$ forms a divided power subalgebra of
$\Cal A_q(n)$, which is $l^n$-dimensional. In this case,
$(x^{(\al)})^l=0$ for any  $x^{(\al)}\in\Cal A_q(n,\bold 1)$ with
$\al\ne 0$  so we say $\Cal A_q(n,\bold 1)$ a {\it quantum
restricted divided power algebra}. \proclaim{Proposition} \ Suppose
that $k$ is a field of characteristic zero, $\ch(q)=l\ge 3$. Then
the quantum divided power algebra $\Cal A_q$ is generated by
elements $x_i$ and $x_i^{(l)}$ $(1\le i\le n)$; and its quantum
restricted divided power algebra $\Cal A_q(n,\bold 1)$ is generated
by $x_i$ $(1\le i\le n)$. In addition, when $l$ is odd, $x_i^{(l)}$
$(1\le i\le n)$ are central elements of $\Cal A_q(n)$ and $\Cal
A_q(n)\cong \Cal A_q(n,\bold 1)\ot_k k[\,x_1^{(l)},\cdots,
x_n^{(l)}\,]$ $($as algebras$)$.
\endproclaim
\demo{Proof}
For any $m\in\Bbb Z_+$, let $m=m_0+m_1l$
such that $0\le M_0\le l-1$, $m_1\ge 0$.
By using Lemma 1.5, we get ${m\bino m_0}=1$, and
conclude from (i) that
$$
\split
\bigl(x_i^{(l)}\bigr)^{m_1}&={2\,l\bino l}{3\,l\bino l}\cdots
{m_1l\bino l}\,x_i^{(m_1l)}=m_1!\,x_i^{(m_1l)},\\
x_i^{(m)}&=x_i^{(m_0)}\,x_i^{(m_1l)}=\frac{1}{m_1!}\,x_i^{(m_0)}\,
\bigl(x_i^{(l)}\bigr)^{m_1}.
\\
x^{(\al)}\,x_i^{(l)}&=\th(\al,l\vn_i)\,x_i^{(l)}\,x^{(\al)}
=x_i^{(l)}\,x^{(\al)} \qquad(\text{\it for \;$q^l=1$}).
\endsplit\tag{\text{ii}}
$$
Note that $\ch(q)=l$ being odd means $q^l=1$.
These imply the statements.
\hfill\qed
\enddemo

By analogy of Corollary 2.3, we have \proclaim {Corollary} The
quantum divided power algebra $\Cal A_q(n)$ in the case when
$\ch(q)=0$ or when $\ch(q)=\text{\it char}\,(k)=p$ $($a prime$)$,
together with the quantum restricted divided power algebra $\Cal
A_q(n,\bold 1)$ in the case when $\ch(q)=\text{\it char}\,(k)=p$, is
a $\Bbb Z^n$-graded braided-commutative Hopf algebra relative to the
bicharacter $\th$ given in {\bf 2.1} $(\text{\rm ii})$. \hfill\qed
\endproclaim

In the latter discussion, we shall prefer to work with the quantum divided
power algebra rather than the quantum $n$-space.

\

\vskip0.35cm
\heading{$\bold 3$. \ $q$-Derivatives, Quantum Groups
and Quantum Weyl Algebra}
\endheading
\vskip0.12cm

\noindent
{\bf 3.1} \ For simplicity of notation, we let $\Cal A_q$ denote
$\Cal A_q(n)$ for any $\ch(q)$, or $\Cal A_q(n,\bold 1)$ only under the case
$\ch(q)=l\,(>2)$.
Consider the algebra automorphisms $\sigma_i$
$(1\le i\le n)$ of $\Cal A_q$ defined by
$$
\sigma_i(x^{(\be)})
=q^{\langle \be,\, \vn_i\rangle }\,x^{(\be)}
=q^{\be_i}\,x^{(\be)}, \qquad \forall \; X^{(\be)}\in\Cal A_q.\leqno(\text{\rm i})
$$
When  $q=1$, one has $\sigma_i=\text{id}$. Obviously, $\sigma_i\sigma_j=
\sigma_j\sigma_i$.

Define $\frac{\pa_q}{\pa\,x_i}$ as the special $q$-derivatives
over $\Cal A_q$ by
$$
\frac{\partial_q}{\partial\,x_i}(x^{(\be)})
=q^{-\vn_i*\be}\,x^{(\be-\vn_i)}, \qquad \forall \; x^{(\be)}\in \Cal A_q.
\leqno(\text{\rm ii})
$$
For convenience,
we briefly let $\pa_i$ denote $\frac{\pa_q}{\pa\,x_i}$.

For $\al\in\Bbb Z^n$, denote $\Th(\al)$ by the algebra automorphisms of
$\Cal A_q$:
$$
\Th(\al)(x^{(\be)})=\th(\al,\be)\,x^{(\be)}, \qquad \forall \; x^{(\be)}\in \Cal AWq,
\leqno(\text{\rm iii})
$$
where $\th$ is the bicharacter on $\Bbb Z^n$ defined in {\bf 2.1}.
Thus we have
\proclaim{Proposition} \ $(1)$ \ $\Th(\al)\Th(\be)=\Th(\al+\be)$,
in particular, $\Th(-\al_i)=\sigma_i\sigma_{i+1}$, for a simple root
$\al_i=\vn_i-\vn_{i+1}$ in root system of type $\bold A_{n-1}$.

\smallskip
$(2)$ \
$\pa_i$ is a $(\Th(-\vn_i)\sigma_i^{\pm 1},\sigma_i^{\mp 1})$-derivation
of $\Cal A_q$, namely,

\smallskip
\hskip.7cm
$\pa_i(x^{(\be)}\,x^{(\ga)})
=\pa_i(x^{(\be)})\,\sigma_i^{\mp 1}(x^{(\ga)})
+\bigl(\Th(-\vn_i)\sigma_i^{\pm 1}\bigr)(x^{(\be)})\,\pa_i(x^{(\ga)}).$

\smallskip
$(3)$ \ $\pa_i\,\pa_j
=\th(-\vn_i,-\vn_j)\,\pa_j\,\pa_i=\th(\vn_i,\vn_j)\,\pa_j\,\pa_i$.

\smallskip
$(4)$ \ $x^{(\al)}(x^{(\be)}\,x^{(\ga)})=(x^{(\al)}\,x^{(\be)})\,x^{(\ga)}$

\smallskip
\hskip2.95cm
$=\th(\al,\be)\,x^{(\be)}\,(x^{(\al)}\,x^{(\ga)})$

\smallskip
\hskip2.95cm
$=\Th(\al)(x^{(\be)})\,(x^{(\al)}\,x^{(\ga)})$.

\smallskip
$(5)$ \ $\sigma_i(x^{(\be)}\,x^{(\ga)})=\sigma_i(x^{(\be)})\,
\sigma_i(x^{(\ga)})$,

\smallskip
\hskip.7cm
$\Th(\al)(x^{(\be)}\,x^{(\ga)})
=\Th(\al)(x^{(\be)})\,\Th(\al)(x^{(\ga)})$.

\smallskip
$(6)$ \
$x^{(\al)}\pa_i$ is a $(\Th(\al-\vn_i)\sigma_i^{\pm 1},\sigma_i^{\mp 1})$-derivation
of $\Cal A_q$, namely,

\smallskip
\noindent
$\bigl(x^{(\al)}\pa_i\bigr)(x^{(\be)}\,x^{(\ga)})
=\bigl(x^{(\al)}\pa_i\bigr)(x^{(\be)})\,\sigma_i^{\mp 1}(x^{(\ga)})
+\bigl(\Th(\al-\vn_i)\sigma_i^{\pm 1}\bigr)(x^{(\be)})\,
\bigl(x^{(\al)}\pa_i\bigr)(x^{(\ga)}).$
\endproclaim
\demo{Proof} \ (1) \ The first claim is due to the property of bicharacter
$\th$ (see {\bf 2.1} (iv)\,). The second follows from $\Th(-\al_i)(x^{(\ga)})=
\th(-\al_i,\ga)x^{(\ga)}=\sigma_i\sigma_{i+1}(x^{(\ga)})$ (by {\bf 2.1} (ii) and
Lemma 2.1 (2) \& (3)\,).

(2) Noting that $[m+m']=[m]\,q^{\mp m'}+[m']\,q^{\pm m}$,
and
$$
{{\be+\ga}\bino{\be}}={{\be+\ga-\vn_i}\bino{\be-\vn_i}}\,q^{\mp\ga_i}+
q^{\pm\be_i}{{\be+\ga-\vn_i}\bino{\be}},
\leqno(\text{\rm iv})
$$
then by {\bf 2.4} (i), {\bf 3.1} (ii) \& Lemma 2.1, we get
$$
\split
\pa_i(&x^{(\be)})\,\sigma_i^{\mp 1}(x^{(\ga)})
+\bigl(\Th(-\vn_i)\sigma_i^{\pm 1}\bigr)(x^{(\be)})\,\pa_i(x^{(\ga)})\\
&=q^{-\vn_i*\be\mp\ga_i}\,x^{(\be-\vn_i)}\,x^{(\ga)}+
\th(-\vn_i,\be)\,q^{\pm\be_i-\vn_i*\ga}x^{(\be)}\,x^{(\ga-\vn_i)}\\
&=\left(\,q^{-\vn_i*\be\mp\ga_i+(\be-\vn_i)*\ga}\,
{{\be+\ga-\vn_i}\bino{\be-\vn_i}}\right.\\
&\qquad\left. + \, q^{-\vn_i*\be+\be*\vn_i\pm\be_i-\vn_i*\ga+\be*(\ga-\vn_i)}\,
{{\be+\ga-\vn_i}\bino{\be}}\,\right)\,x^{(\be+\ga-\vn_i)}\\
&=q^{-\vn_i*(\be+\ga)+\be*\ga}\,
\left({{\be+\ga-\vn_i}\bino{\be-\vn_i}}\,q^{\mp\ga_i}+
q^{\pm\be_i}\,{{\be+\ga-\vn_i}\bino{\be}}\,\right)\,x^{(\be+\ga-\vn_i)}\\
&=q^{-\vn_i*(\be+\ga)+\be*\ga}\,{{\be+\ga}\bino{\be}}\,x^{(\be+\ga-\vn_i)}\\
&=q^{\be*\ga}\,{{\be+\ga}\bino{\be}}\,\pa_i(x^{(\be+\ga)})
=\pa_i(x^{(\be)}\,x^{(\ga)}).
\endsplit
$$
Therefore, $\pa_i$ is a
$(\Th(-\vn_i)\sigma_i^{\pm 1},\sigma_i^{\mp 1})$-derivation of
$\Cal A_q$ in the sense of {\bf 1.3} (iii).

(3) By (ii) and {\bf 2.1} (ii), for any $x^{(\ga)}\in\Cal A_q$, we get
$$
\split
\pa_i\,\pa_j(x^{(\ga)})
&=q^{-(\vn_i+\vn_j)*\ga+\vn_i*\vn_j}\,x^{(\ga-\vn_i-\vn_j)}\\
&=\th(\vn_i,\vn_j)\,q^{-(\vn_i+\vn_j)*\ga+
\vn_j*\vn_i}\,x^{(\ga-\vn_i-\vn_j)}\\
&=\th(\vn_i,\vn_j)\,\pa_j\,\pa_i(x^{(\ga)}).
\endsplit
$$

(4) follows from the associativity and the $\th$-commutativity
of $\Cal A_q$.

(5) is clear.

(6) is obtained by combining (2) with (4), and observing the additivity
of $\Th$ in (1).
\hfill
\qed
\enddemo

\noindent
{\it Remark.} \ The commutative coefficients between $\pa_i$ in
Proposition 3.1 (3) coincide with those of $x^{(\vn_i)}=x_i$, that is,
$\pa_i\,\pa_j=\th(\vn_i,\vn_j)\,\pa_j\,\pa_i$ and $x_i\,x_j=\th(\vn_i,\vn_j)
x_j\,x_i$.
In the discussion of Wess \&
Zumino in \cite{17}, however, $x\,y=q\,y\,x$ but $\pa_x\,\pa_y=q^{-1}\pa_y\,\pa_x$
(cf. p. 309 \cite{17}, (4.5) \& (4.7)), which distinguishes from our case.
Especially, in process of our introducing the $q$-differential operators
$\pa_i$ over $\Cal A_q$, of particular interest
is that will lead to a new quantum group structure below (in the sense
of Drinfeld).

\medskip
\noindent
{\bf 3.2} \
Now let $\frak D_q$ be the associative algebra over $k$
generated by
the symbols $\Th(\pm\vn_i)$, $\sigma_i^{\pm 1}$, $\pa_i$ ($1\le i\le n$),
associated to the bicharacter $\th$ on $\Bbb Z^n$ given in {\bf 2.1},
satisfying the following relations:
$$
\gather
\sigma_i\sigma_i^{-1}=1=\sigma_i^{-1}\sigma_i, \qquad
\sigma_i\sigma_j=\sigma_j\sigma_i,\tag{i}\\
\sigma_i\,\Th(\vn_j)=\Th(\vn_j)\,\sigma_i, \qquad
\Th(-\vn_i+\vn_{i+1})=\sigma_i\sigma_{i+1},\\
\Th(\vn_i)^{-1}=\Th(-\vn_i), \qquad \Th(0)=1,\tag{ii}\\
\Th(\vn_i)\,\Th(\vn_j)=\Th(\vn_i+\vn_j)=\Th(\vn_j)\,\Th(\vn_i),\\
\Th(\vn_j)\,\pa_i\,\Th(\vn_j)^{-1}=\th(\vn_i,\vn_j)\,\pa_i,\tag{iii}\\
\sigma_j\,\pa_i\,\sigma_j^{-1}=q^{-\delta_{ij}}\pa_i,\tag{iv}\\
\pa_i\,\pa_j=\th(\vn_i,\vn_j)\,\pa_j\,\pa_i.\tag{v}
\endgather
$$

Furthermore, $\frak D_q$ can be equipped with a quantum group structure
if we define the following mappings $\Delta$, $\ep$ and $S$
on the generators of $\frak D_q$ as
$$
\gather
\Delta: \frak D_q\to \frak D_q\ot\frak D_q\tag{vi}\\
\Delta\,(\sigma_i^{\pm 1})=\sigma_i^{\pm 1}\ot \sigma_i^{\pm 1}, \\
\Delta\,(\Th(\pm\vn_i))=\Th(\pm\vn_i)\ot\Th(\pm\vn_i),\\
\Delta\,(\pa_i)=\pa_i\ot\sigma_i^{-1}+\Th(-\vn_i)\sigma_i\ot\pa_i.\\
\ep: \frak D_q\to k\tag{vii}\\
\ep(\sigma_i^{\pm 1})=1=\ep(\Th(\pm\vn_i)),\\
\ep(\pa_i)=0.\\
S: \frak D_q\to \frak D_q\tag{viii}\\
S(\sigma_i^{\pm 1})=\sigma_i^{\mp 1}, \\
S(\Th(\pm\vn_i))=\Th(\mp\vn_i),\\
S(\pa_i)=-q\,\Th(\vn_i)\,\pa_i.
\endgather
$$

Again we extend the definitions of $\Delta$, $\ep$ (resp. $S$) on
$\frak D_q$ (anti-)algebraically. Thus we obtain the following
\proclaim{Theorem} \ $(\frak D_q, \Delta, \ep, S)$ is a quantum
group with the above relations  $(\text{\rm i})$ --- $(\text{\rm
viii})$.
\endproclaim
\demo{Proof} First we need to show that $\Delta$, $\ep$ and $S$
preserve the algebraic relations (i)---(v) of $\frak D_q$. This is clear
for $\ep$ and $S$, and clear for $\Delta$ preserving relations (i)---(ii).

So it remains to check it for $\Delta$ with respect to
relations (iii)---(v). Note that
$$
\gather
(\Th(-\vn_i)\si_i)\,\pa_j=\th(\vn_i,\vn_j)\,
\pa_j\,(\Th(-\vn_i)\si_i),\\
\pa_i\,(\Th(-\vn_j)\si_j)=\th(\vn_i,\vn_j)\,(\Th(-\vn_j)\si_j)\,\pa_i,
\endgather
$$
for $i\ne j$. Hence, we have
$$
\split
\Delta(\Th(\vn_j))\Delta(\pa_i)\Delta(\Th(\vn_j)^{-1})
&=\Th(\vn_j)\,\pa_i\,\Th(\vn_j)^{-1}\ot\sigma_i^{-1}\\
&\qquad + \Th(-\vn_i)\sigma_i\ot
\Th(\vn_j)\,\pa_i\,\Th(\vn_j)^{-1}\\
&=\th(\vn_i,\vn_j)\Delta(\pa_i),\endsplit$$
$$\split
\Delta(\sigma_j)\Delta(\pa_i)\Delta(\sigma_j^{-1})
&=\si_j\,\pa_i\,\si_j^{-1}\ot \si_i^{-1}+\Th(-\vn_i)\si_i\ot\si_j\,\pa_i\,\si_j^{-1}\\
&=q^{-\delta_{ij}}\Delta(\pa_i),\endsplit$$
$$\split
\De(\pa_i)\De(\pa_j)
&=(\pa_i\ot \si_i^{-1}+\Th(-\vn_i)\si_i\ot\pa_i)
(\pa_j\ot \si_j^{-1}+\Th(-\vn_j)\si_j\ot\pa_j)\\
&=\pa_i\pa_j\ot\si_i^{-1}\si_j^{-1}+\Th(-\vn_i-\vn_j)\si_i\si_j\ot\pa_i\pa_j\\
&\qquad + \,
(\Th(-\vn_i)\si_i)\,\pa_j\ot\pa_i\si_j^{-1}+\pa_i\,(\Th(-\vn_j)\si_j)\ot\si_i^{-1}\pa_j\\
&=\th(\vn_i,\vn_j)\bigl(\pa_j\pa_i\ot\si_j^{-1}\si_i^{-1}
+\Th(-\vn_i-\vn_j)\si_j\si_i\ot \pa_j\pa_i\bigr.\\
&\qquad +\,
\bigl.(\Th(-\vn_j)\si_j)\,\pa_i\ot\pa_j\si_i^{-1}
+\pa_j\,(\Th(-\vn_i)\si_i)\ot\si_j^{-1}\pa_i\bigr)\\
&=\th(\vn_i,\vn_j)\De(\pa_j)\De(\pa_i), \qquad (i\ne j).
\endsplit
$$

In view of the fact just proved, together with (vi) \& (vii), we see that
$(1\ot\Delta)\Delta=(\Delta\ot 1)\Delta$ and
$(1\ot\ep)\Delta=1=(\ep\ot 1)\Delta$ hold. Again by (vi) \& (viii), we
have
$m\circ(1\ot S)\circ\De(\pa_i)
=\pa_i\si_i+\Th(-\vn_i)\,\si_i\,(-q\Th(\vn_i)\pa_i)
=\pa_i\si_i-q\si_i\pa_i=0$ and
$m\circ(S\ot 1)\circ\De(\pa_i)=-q\Th(\vn_i)\,\pa_i\,\si_i^{-1}
+\si_i^{-1}\,\Th(\vn_i)\,\pa_i
=\Th(\vn_i)(-q\pa_i\si_i^{-1}+\si_i^{-1}\pa_i)=0$,
thus we get $m\circ(1\ot S)\circ\De(\pa_i)=\ep(\pa_i)
=m\circ(S\ot 1)\circ\De(\pa_i)$.
On the other hand, owing to (vi) and
$\Th(\pm\vn_i)\,\Th(\mp\vn_i)=1=\si_i^{\pm 1}\si_i^{\mp 1}$,
there holds
$m\circ(1\ot S)\circ\De=\eta\circ\ep=m\circ(S\ot 1)\circ\De$,
where $(\frak D_q, m, \eta)$ is the algebra structure of $\frak D_q$.

Thereby, $(\frak D_q, m, \eta, \De, \ep, S)$ is a non-commutative and
non-cocommutative Hopf algebra, namely, a quantum group.
\hfill\qed
\enddemo

\noindent
{\it Remark.} \
Actually, we can equip $\frak D_q$ with another quantum group structure
$(\frak D_q$, $\Delta^{(-)}$, $\ep$, $S^{(-)}$), where only one difference
is the actions
of $\De^{(-)}$ and $S^{(-)}$ on $\pa_i$ ($1\le i\le n$)
given respectively by $\De^{(-)}(\pa_i)
=\pa_i\ot\si_i+\Th(-\vn_i)\si_i^{-1}\ot\pa_i$, and
$S^{(-)}(\pa_i)=-q^{-1}\,\Th(\vn_i)\,\pa_i$.

\medskip
\noindent
{\bf 3.3} \ Based on the structure of quantum group $\frak D_q$
in Theorem 3.2 and observing that the commutative rule
for $x_i$ in $\Cal A_q$ is the same as that of the special
$q$-derivatives $\pa_i$ (see Remark 3.1),
we can augment $\Cal A_q$ through adding a certain multiplication abelian
group $\Theta$ (as group-like elements), and construct another
quantum group $\frak A_q$ below such that it contains the quantum
divided power algebra $\Cal A_q$ as its subalgebra.
Here
$\Theta=\{\Th(\al)\mid \al\in\Bbb Z^n\}$ acts
conjugately on $\Cal A_q$ as its an automorphism group.
When $\Cal A_q$ (as an object in the braided category $\Cal {GB}$)
is of a braided Hopf algebra structure, we has a
reasonable interpretation for such a construction of the quantum group
$\frak A_q$ (as an object in the usual category $\Cal {HA}$ of Hopf
algebras): the price introducing
the group $\Theta$ into $\Cal A_q$ lies in transmuting the
``braided'' twisting of the algebra structure on $\Cal A_q\otimes \Cal A_q$
(in $\Cal {GB}$) into the ``trivial'' twisting of the algebra structure
on $\frak A_q\otimes \frak A_q$ (in $\Cal {HA}$). As is clear, according to
their respective comultiplications (as algebra homomorphisms!).

Let $\frak A_q$ be the associative algebra over $k$ generated by the
quantum divided power (restricted) algebra $\Cal A_q$, together with the symbols
$\Th(\pm\vn_i)$ ($1\le i\le n$), associated to the bicharacter $\th$ on
$\Bbb Z^n$ given in {\bf 2.1}, subject to the relations
$$
\gather
\Th(\vn_i)^{-1}=\Th(-\vn_i), \qquad \Th(0)=1,\tag{i}\\
\Th(\vn_i)\,\Th(\vn_j)=\Th(\vn_i+\vn_j)=\Th(\vn_j)\,\Th(\vn_i),\\
\Th(\vn_j)\,x_i\,\Th(\vn_j)^{-1}=\th(\vn_j,\vn_i)\,x_i,\tag{ii}\\
\bigl( \,\Th(\vn_j)\,x^{(l)}_i\,\Th(\vn_j)^{-1}=x_i^{(l)}\quad
\text{\it in $\frak  A_q(n)$ if in addition $q^l=1$} \,\bigr)\\
x_i\,x_j=\th(\vn_i,\vn_j)\,x_j\,x_i.\tag{iii}\\
\bigl( \,x_i\,x^{(l)}_j=x^{(l)}_j\,x_i\quad
\text{\it in $\frak  A_q(n)$ if in addition $q^l=1$} \,\bigr)
\endgather
$$
Note that $\frak A_q$ here also indicates $\frak A_q(n)$ for any $\ch(q)$
or $\frak A_q(n,\bold 1)$ only for $q^l=1$.
Similar to Theorem 3.2, we have
\proclaim
{Theorem} \
$(\frak A_q, \De, \ep, S)$ forms a quantum group, which is
the required quantization object of the polynomial algebra in $n$
variables in the context of quantum groups when $\ch(q)=0$ $($also can be
viewed as the quantized
universal enveloping algebra of the abelian Lie algebra of dimension $n$\,$)$,
under the comultiplication $\De$, the counit $\ep$ and the antipode $S$
given by
$$
\gather
\Delta: \frak A_q\to \frak A_q\ot\frak A_q\tag{iv}\\
\Delta\,(\Th(\pm\vn_i))=\Th(\pm\vn_i)\ot\Th(\pm\vn_i),\\
\Delta\,(x_i)=x_i\ot 1+\Th(\vn_i)\ot x_i.\\
\bigl( \,\De\,(x^{(l)}_i)=x^{(l)}_i\ot 1+1\ot x^{(l)}_i\quad
\text{\it in $\frak A_q(n)$ only for $q^l=1$} \,\bigr)\\
\ep: \frak A_q\to k\tag{v}\\
\ep(\Th(\pm\vn_i))=1,\\
\ep(x_i)=0. \\
\bigl( \,\ep(x_i^{(l)})=0\quad
\text{\it in $\frak A_q(n)$ only for $q^l=1$} \,\bigr)\\
S: \frak A_q\to \frak A_q\tag{vi}\\
S(\Th(\pm\vn_i))=\Th(\mp\vn_i),\\
S(x_i)=-\Th(-\vn_i)\,x_i.\\
\bigl( \,S(x_i^{(l)})=-x_i^{(l)}\quad
\text{\it in $\frak A_q(n)$ only for $q^l=1$} \,\bigr)
\endgather
$$
\endproclaim
\demo{Proof} It follows from a similar argument as Theorem 3.2.
\hfill\qed
\enddemo

\noindent
{\it Remark.} \
Given a Hopf algebra $(H, \De, \ep, S)$, for $a\in H$,
consider its {\it left adjoint action}
$\text{ad}\,a$ on $H$: $\text{ad}\,a(b)=a_1bS(a_2)$ where
$\De(a)=a_1\ot a_2$ (in Sweedler convention).
Clearly, the map $a\mapsto \text{ad}\,a$ of $H$ into $\text{End}\,H$ is an
algebra homomorphism. Now for $H=\frak A_q$, consider the action of
$\text{ad}\,x_i$ on its subalgebra $\Cal A_q$, that is, $\forall\,
x^{(\ga)}\in\Cal A_q$, we have
$\text{ad}\,x_i(x^{(\ga)})
=x_i\,x^{(\ga)}-\Th(\vn_i)\,x^{(\ga)}\,\Th(-\vn_i)\,x_i
=x_i\,x^{(\ga)}-\th(\vn_i,\ga)\,x^{(\ga)}\,x_i
=0$.
On the other hand, if in addition $q^l=1$, we get
$\text{ad}\,x^{(l)}_i(x^{(\ga)})=x^{(l)}_i\,x^{(\ga)}-x^{(\ga)}\,x^{(l)}_i=0$
in $\Cal A_q(n)$.
So $\text{ad}\,x^{(\al)}|_{\Cal A_q}=0$ for any $x^{(\al)}\in
\Cal A_q$. This fact is compatible with the classical case, since
the group $\Theta$ degenerates into the unit group and $\frak A_q(n)$
into a polynomial algebra in $n$ variables with a known Hopf algebra
structure when $q$ takes $1$. Consequently, in the case when
$\ch(q)=0$, the quantum group $\frak A_q(n)$
achieved above is just the corresponding {\it object} of the
polynomial algebra in $n$ variables in the context of quantum groups,
which also can be considered as the quantized universal
enveloping algebra of the abelian Lie algebra of dimension $n$.

\medskip
\noindent
{\bf 3.4} \ By analogy of the argument in {\bf 3.2},
we denote $\frak U_q$ by
the associative algebra over $k$ generated by the symbols
$\Th(\pm\vn_i)$, $\sigma_i^{\pm 1}$, $x_i$ ($1\le i\le n$),
associated to the bicharacter $\th$ on $\Bbb Z^n$ given in {\bf 2.1},
satisfying the following relations:
$$
\gather
\sigma_i\sigma_i^{-1}=1=\sigma_i^{-1}\sigma_i, \qquad
\sigma_i\sigma_j=\sigma_j\sigma_i,\tag{i}\\
\sigma_i\,\Th(\vn_j)=\Th(\vn_j)\,\sigma_i, \qquad
\Th(-\vn_i+\vn_{i+1})=\sigma_i\sigma_{i+1},\\
\Th(\vn_i)^{-1}=\Th(-\vn_i), \qquad \Th(0)=1,\tag{ii}\\
\Th(\vn_i)\,\Th(\vn_j)=\Th(\vn_i+\vn_j)=\Th(\vn_j)\,\Th(\vn_i),\\
\bigl(\,\Th(\vn_i)^l=\Th(l\,\vn_i)=1, \qquad \sigma_i^l=1 \quad \text{\it
only when} \quad q^l=1\,\bigr)\\
\Th(\vn_i)\,x_j\,\Th(\vn_i)^{-1}=\th(\vn_i,\vn_j)\,x_j,\tag{iii}\\
\sigma_i\,x_j\,\sigma_i^{-1}=q^{\delta_{ij}}\,x_j,\tag{iv}\\
x_i\,x_j=\th(\vn_i,\vn_j)\,x_j\,x_i.\tag{v}
\endgather
$$

Moreover, the comultiplication $\De$, the counity $\ep$ and the
antipode $S$ over $\frak U_q$ are defined respectively by
$$
\gather
\Delta: \frak U_q\to \frak U_q\ot\frak U_q\tag{vi}\\
\Delta\,(\sigma_i^{\pm 1})=\sigma_i^{\pm 1}\ot \sigma_i^{\pm 1}, \\
\Delta\,(\Th(\pm\vn_i))=\Th(\pm\vn_i)\ot\Th(\pm\vn_i),\\
\Delta\,(x_i)=x_i\ot\sigma_i+\Th(\vn_i)\sigma_i^{-1}\ot x_i.\\
\ep: \frak U_q\to k\tag{vii}\\
\ep(\sigma_i^{\pm 1})=1=\ep(\Th(\pm\vn_i)),\\
\ep(x_i)=0.\\
S: \frak U_q\to \frak U_q\tag{viii}\\
S(\sigma_i^{\pm 1})=\sigma_i^{\mp 1}, \\
S(\Th(\pm\vn_i))=\Th(\mp\vn_i),\\
S(x_i)=-q\,\Th(-\vn_i)\,x_i.\\
\endgather
$$

Similar to Theorem 3.2, we obtain \proclaim{Theorem} \ $(\frak U_q,
\Delta, \ep, S)$ is the quantum group corresponding to the quantum
$n$-space $k[A_q^{n|0}]$  in particular, in the case when $\ch(q)=l$
is odd, whose Hopf algebra structure restricted on its central
subalgebra $k[x_1^l,\cdots,x_n^l]$ is just the usual Hopf algebra
structure of the polynomial algebra in $n$ variables.
\endproclaim
\demo{Proof} \
The first claim follows from a similar argument as Theorem 3.2.

On the other hand,
using formula (vi) and induction on $m$, we get
$$
\De(x_i^m)=\sum_{k=0}^m{m\bino k}\,x_i^{m-k}\,\Th(\vn_i)^k\,\si_i^{-k}\ot
x_i^k\,\si_i^{m-k}.
$$
In particular, since $\ch(q)=l$ being odd means $q^l=1$, we have
${l\bino k}=0$ for $1\le k<l$ and
$\Th(\vn_i)^l\,\si_i^{-l}=1$ so that $\De(x_i^l)=x_i^l\ot 1+1\ot x_i^l$
and $S(x_i^l)=-x_i^l$.
These are compatible with the fact that $k[x_1^l,\cdots,x_n^l]$ is as
a (polynomial) central subalgebra of $k[A_q^{n|0}]$ when $q^l=1$.
\hfill\qed
\enddemo

\noindent
{\it Remark.} \ Actually, since the fact that $x_i \,(1\le i\le n)$ are the
generators of $k[A_q^{n|0}]$ is independently of the {\it characteristic}
$\ch(q)$ of $q$. Thereby, when $q$ is a primitive $l$-th root of unity,
the quantum group object corresponding to
the quantum $n$-space $k[A_q^{n|0}]$ can be directly defined by the same
relations (i) --- (viii) as in {\bf 3.4}.

\medskip
\noindent
{\bf 3.5} \ Let $\frak D_q^{(\pm)}:=(\frak D_q, \De^{(\pm)}, \ep,
S^{(\pm)})$ denote the two Hopf algebras appearing in {\bf 3.2}. Then from
Proposition 3.1 (2) \& (5), we see that $\Cal A_q$ is a (left)
$\frak D_q^{(\pm)}$-module algebra. So we are able to make their smash
product algebras $\Cal A_q\#\frak D_q^{(\pm)}$ in a familiar fashion as in
\cite{16}, which are the same as
$\Cal A_q\ot\frak D_q$ as vector spaces but multiplications respectively
given by
$$
\split
(x^{(\al)}\,\# \,\pa_i)\circ(x^{(\be)}\,\#\, d)
&=x^{(\al)}\,\pa_i(x^{(\be)})\,\#\,\si_i^{\mp 1}\, d
+x^{(\al)}(\Th(-\vn_i)\si_i^{\pm 1})(x^{(\be)})\#\,\pa_i\, d\\
&=\,q^{\al*\be-\al*\vn_i-\vn_i*\be}\,{{\al+\be-\vn_i}\bino{\al}}
x^{(\al+\be-\vn_i)}\# \si_i^{\mp 1}\, d\\
&\qquad + \;\th(\be,\vn_i)\,q^{\al*\be\pm\be_i}\,{{\al+\be}\bino{\al}}
x^{(\al+\be)}\,\#\, \pa_i\, d, \\
(x^{(\al)}\,\#\, g)\circ(x^{(\be)}\,\#\, d)
&=x^{(\al)}\,g(x^{(\be)})\,\#\, g\,d,
\endsplit
$$
where
$\De^{(\pm)}(\pa_i)=\pa_i\ot\si_i^{\mp 1}+\Th(-\vn_i)\,\si_i^{\pm 1}\ot\pa_i$ and
$\De^{(\pm)}(g)=g\ot g$,
for $x^{(\al)}, x^{(\be)}\in\Cal A_q$ and $\pa_i, g, d\in\frak D_q$.
More precisely, we have
$$
\split
(x^{(\al)}\,\#\, 1)\circ(1\,\#\,d)&=x^{(\al)}\,\#\,d,\\
(1\,\#\,\Th(\vn_i))\circ(x_j\,\#\,1)&=\Th(\vn_i)(x_j)\,\#\,\Th(\vn_i)\\
&=\th(\vn_i,\vn_j)\,x_j\,\#\,\Th(\vn_i),\\
(1\,\#\,\si_i)\circ(x_j\,\#\,1)&=q^{\de_{ij}}\,x_j\,\#\,\si_i,\\
(1\,\#\,\pa_i)\circ(x_j\,\#\,1)&=\de_{ij}\,\#\,\si_i^{\mp 1}+
\th(\vn_j,\vn_i)\,q^{\pm\de_{ij}}\,x_j\,\#\,\pa_i.
\endsplit
\tag{i}
$$

Actually, if we briefly identify elements $x^{(\al)}\# d$ in
$\Cal A_q\#\frak D_q^{(\pm)}$ with $x^{(\al)}\,d$, then the smash product
algebras $\Cal A_q\#\frak D_q^{(\pm)}$ containing $\Cal A_q$ and
$\frak D_q^{(\pm)}$ as subalgebras are just the {\it quantum
differential operators algebras} over $\Cal A_q$, which will degenerate
into the usual differential operators algebra when  $q$ takes $1$.
Particularly, combining formulae (i) with the identification above,
we get the relations below
$$
\split
x^{(\al)}\circ d&=x^{(\al)}\,d,\\
\Th(\vn_i)\circ x_j\circ\Th(\vn_i)^{-1}&=\th(\vn_i,\vn_j)\,x_j,\\
\si_i\circ x_j\circ\si_i^{-1}&=q^{\de_{ij}}\,x_j,\\
\pa_i\circ x_j&=\de_{ij}\si_i^{\mp 1}
+\th(\vn_j,\vn_i)\,q^{\pm\de_{ij}}x_j\circ\pa_i.
\endsplit\tag{\text{\rm ii}}
$$

Applying formulae (ii), together with $\Cal A_q$ and $\frak D_q^{(\pm)}$,
we are able to
construct the following {\it quantum Weyl algebra},
which is different from those appeared in the literature  for instance,
Proposition 5.2.2 in \cite{3}, Section 2.1 in \cite{5}, \cite{10}, etc.).
\proclaim{Definition}
Let $\Cal W_q(2n)$ be the associative algebra over $k$ generated by
the symbols $\Th(\pm\vn_i)$, $\si_i^{\pm 1}$, $x_i$ and $\pa_i$
$(1\le i\le n)$, associated to the bicharacter $\th$ on $\Bbb Z^n$
defined in {\bf 2.1}, obeying the following relations:
$$
\gather
\Th(\pm\vn_i)\circ\Th(\mp\vn_i)=1=\si_i^{\pm 1}\circ\si_i^{\mp 1},\tag{iii} \\
\Th(-\vn_i+\vn_{i+1})=\si_i\circ\si_{i+1},\\
\Th(\vn_i)\circ\Th(\vn_j)=\Th(\vn_i+\vn_j)=\Th(\vn_j)\circ\Th(\vn_i), \\
\si_i\circ\si_j=\si_j\circ\si_i, \\
\si_i\circ\Th(\vn_j)=\Th(\vn_j)\circ\si_i,\\
\Th(\vn_i)\circ x_j\circ\Th(-\vn_i)=\th(\vn_i,\vn_j)\,x_j,\tag{iv} \\
\Th(\vn_i)\circ\pa_j\circ\Th(-\vn_i)=\th(\vn_j,\vn_i)\,\pa_j,\\
\si_i\circ x_j\circ\si_i^{-1}=q^{\de_{ij}}\,x_j,\tag{v} \\
\si_i\circ\pa_j\circ\si_i^{-1}=q^{-\de_{ij}}\,\pa_j,\\
x_i\circ x_j
=\th(\vn_i,\vn_j)\,x_j\circ X_i,\tag{vi} \\
\pa_i\circ\pa_j=\th(\vn_i,\vn_j)\,\pa_j\circ\pa_i,\\
\pa_i\circ x_j=\th(\vn_j,\vn_i)\,x_j\circ \pa_i,
\qquad (i\ne j),\tag{vii}\\
\pa_i\circ X_i-q^{\pm 1}x_i\circ\pa_i=\si_i^{\mp 1}.\tag{viii}\\
\endgather
$$
where the relations $(\text{\rm viii})$ are equivalent to the following relations:
$$
\pa_i\circ x_i=\frac{q\,\si_i-(q\,\si_i)^{-1}}{q-q^{-1}}, \qquad
x_i\circ \pa_i=\frac{\si_i-\si_i^{-1}}{q-q^{-1}}.\leqno(\text{\rm ix})
$$
\endproclaim

\noindent
{\it Remark.} \ The relations above imply $\Cal W_q(2n)$
also contains
$\frak U_q$ given in {\bf 3.4} as its subalgebra. However, $\Cal W_q(2n)$
itself doesn't to be a Hopf algebra due to the last relations (viii),
which is coincident with the classical situation.

\

\medskip
\bigskip
\heading{$\bold 4$. \ An Application: Realization}
\endheading

\medskip
\noindent
{\bf 4.1} \
As an application of discussions in section {\bf 3}, we are now in the
position to realize the quantized
universal enveloping algebra $U_q(\frak g)$
(where $\frak g=\frak {gl}_n$ or $\frak {sl}_n$)
as certain $q$-differential operators in $\Cal W_q(2n)$
defined over the quantum divided power (restricted) algebra $\Cal A_q$,
such that the quantum divided power (restricted) algebra $\Cal A_q$ is
made into a $U_q(\frak g)$-module algebra in the sense of {\bf 1.3}.

\proclaim{Theorem} \ For any monomial $x^{(\be)}\in \Cal A_q$ and
$1\le i<n$, set
$$
\gather
e_i(x^{(\be)})=\bigl(\,x_i\,\pa_{i+1}\,\sigma_i\,\bigr)\,(x^{(\be)}),\tag{i}\\
f_i(x^{(\be)})=\bigl(\,\sigma_i^{-1}\,x_{i+1}\,\pa_i\,\bigr)\,(x^{(\be)}),\tag{ii}\\
\Cal K_i(x^{(\be)})=(\sigma_i\,\sigma_{i+1}^{-1})\,(x^{(\be)}),\tag{iii}\\
\Cal K_i^{-1}(x^{(\be)})=(\sigma_i^{-1}\,\sigma_{i+1})\,(x^{(\be)}).\tag{iv}
\endgather
$$
Formulas $(\text{\rm i})$--$(\text{\rm iv})$ define the structure of a $U_q(\frak {sl}_n)$-module
algebra on $\Cal A_q$.
\qed
\endproclaim
\demo{Proof} \ The proof will be given in two steps.

(I) \ We first show that the formulas (i)--(iv) equip $\Cal A_q$
with a $U_q(\frak {sl}_n)$-module structure. To do this, we need to
check the algebra relations {\bf 1.4} (i)--(v) of $U_q(\frak
{sl}_n)$.

Using Lemma 2.1 (2), {\bf 2.4} (i) \& {\bf 3.1} (ii), we get from
relations (i) \& (ii) that
$$
\split
e_i(x^{(\be)})&=q^{\be_i-\vn_{i+1}*\be+\vn_i*(\be-\vn_{i+1})}\,
{\be+\vn_i-\vn_{i+1}\bino \vn_i}\,x^{(\be+\vn_i-\vn_{i+1})}\\
&=[\,\be_i+1\,]\,x^{(\be+\vn_i-\vn_{i+1})},\\
f_i(x^{(\be)})&=q^{-\vn_i*\be+\vn_{i+1}*(\be-\vn_i)-(\be_i-1)}\,
{\be-\vn_i+\vn_{i+1}\bino \vn_{i+1}}\,x^{(\be-\vn_i+\vn_{i+1})}\\
&=[\,\be_{i+1}+1\,]\,x^{(\be-\vn_i+\vn_{i+1})}.
\endsplit\tag{v}
$$
Clearly, relation {\bf 1.4} (i) holds.
For relation {\bf 1.4} (ii), we
have
$$
\split
\Cal K_i\,e_j\,\Cal K_i^{-1}(x^{(\be)})
&=q^{-\be_i+\be_{i+1}}\,[\be_j+1]\,\Cal K_i\,(x^{(\be+\vn_j-\vn_{j+1})})\\
&=q^{\delta_{ij}-\delta_{i,j+1}-\delta_{i+1,j}+\delta_{i+1,j+1}}\,
  e_j(x^{(\be)}).
\endsplit
$$
We can show $\Cal K_i\,f_j\,\Cal K_i^{-1}=q^{-a_{ij}}f_j$ in a similar fashion.
For relation {\bf 1.4} (iii) we have
$$
\split
[e_i,f_j]\,(x^{(\be)})
&=[\be_{j+1}+1]\;e_i(x^{(\be-\vn_j+\vn_{j+1})})
 -[\be_i+1]\;f_j(x^{(\be+\vn_i-\vn_{i+1})})\\
&=\bigl(\;[\be_{j+1}+1]\;[\be_i+1-\delta_{j,i}+\delta_{j+1,i}]\\
&\qquad -[\be_i+1]\;[\be_{j+1}+1+\delta_{i,j+1}-\delta_{i+1,j+1}]\;\bigr)\;
 x^{(\be+\vn_i+\vn_{j+1}-\vn_{i+1}-\vn_j)}\\
&=\delta_{ij}\bigl(\;[\be_{i+1}+1]\;[\be_i]-[\be_i+1]\;[\be_{i+1}]\;\bigr)\;x^{(\be)}\\
&=\delta_{ij}\frac{q^{\be_i-\be_{i+1}}-q^{\be_{i+1}-\be_i}}{q-q^{-1}}\;x^{(\be)}\\
&=\delta_{ij}\frac{\Cal K_i-\Cal K_i^{-1}}{q-q^{-1}}\;(x^{(\be)}).
\endsplit
$$
As for relation {\bf 1.4} (iv), if $|i-j|>1$, then from
formula (v) of $e_i$ it is easy to see that $e_ie_j=e_je_i$. If $|i-j|=1$,
without loss of generality, assume that $j=i+1$, then we shall check that
$$
\bigl(\;e_i^2\,e_{i+1}-(q+q^{-1})\,e_i\,e_{i+1}\,e_i+e_{i+1}\,e_i^2\;\bigr)\;(x^{(\be)})=0.
\leqno(*)
$$
Observing that
$$
\split
e_i^2\;e_{i+1}(x^{(\be)})&=[\be_{i+1}+1]\;[\be_i+1]\;[\be_i+2]\,
x^{(\be+2\vn_i-\vn_{i+1}-\vn_{i+2})},\\
e_i\;e_{i+1}\;e_i(x^{(\be)})&=[\be_i+1]\;[\be_{i+1}]\;[\be_i+2]\,
x^{(\be+2\vn_i-\vn_{i+1}-\vn_{i+2})},\\
e_{i+1}\;e_i^2(x^{(\be)})&=[\be_i+1]\;[\be_i+2]\;[\be_{i+1}-1]\,
x^{(\be+2\vn_i-\vn_{i+1}-\vn_{i+2})},\\
\endsplit
$$
and $[m+1]-(q+q^{-1})\;[m]+[m-1]=0$, we see that the equality (*) holds.
Similarly, we can prove
$$
\bigl(\;e_{i+1}^2\;e_i-(q+q^{-1})\;e_{i+1}\;e_i\;e_{i+1}+e_i\;e_{i+1}^2\;\bigr)\;(x^{(\be)})=0.
$$
As for the last relation {\bf 1.4} (v), it is clear from formula (v)
that $f_if_j=f_jf_i$ if $|i-j|>1$. When $|i-j|=1$, by (v), we have
$$\split
\bigl(\;f_{i+1}^2\;f_i&-(q+q^{-1})\;f_{i+1}\;f_i\;f_{i+1}+f_i\;f_{i+1}^2\;\bigr)\;(x^{(\be)})\\
&=\bigl(\;[\be_{i+1}+1]\;[\be_{i+2}+1]\;[\be_{i+2}+2]\bigr.\\
&\qquad -\,(q+q^{-1})\,[\be_{i+2}+1]\;[\be_{i+1}]\;[\be_{i+2}+2]\\
&\qquad +\,\bigl. [\be_{i+2}+1]\;[\be_{i+2}+2]\;[\be_{i+1}-1]\;\bigr)\,
x^{(\be-\vn_i-\vn_{i+1}+2\vn_{i+2})}=0.
\endsplit
$$
Similarly, we can check that
$\bigl(\;f_i^2\;f_{i+1}-(q+q^{-1})\;f_i\;f_{i+1}\;f_i
+f_{i+1}\;f_i^2\;\bigr)\;(x^{(\be)})=0$.

\medskip
(II) \ We next prove that the quantum divided power (restricted)
algebra $\Cal A_q$ is a $U_q(\frak {sl}_n)$-algebra.
By {\bf 1.3} (i) \& (ii), {\bf 1.4} (vi)--(viii) and noting
that the definition of $\Cal K_i^{\pm}$ in (iii)--(iv), we need
to check that for any $u\in U_q(\frak {sl}_n)$, there hold
$$
\gather
u\,1=\epsilon(u)\;1,\tag{1}\\
\Cal K_i(x^{(\be)}\;x^{(\ga)})
=\Cal K_i(x^{(\be)})\;\Cal K_i(x^{(\ga)}),\tag{2}\\
e_i(x^{(\be)}\;x^{(\ga)})
=x^{(\be)}\;e_i(x^{(\ga)})+e_i(x^{(\be)})\;\Cal K_i(x^{(\ga)}),\tag{3}\\
f_i(x^{(\be)}\;x^{(\ga)})
=\Cal K_i^{-1}(x^{(\be)})\;f_i(x^{(\ga)})+f_i(x^{(\be)})\;x^{(\ga)},
\tag{4}
\endgather
$$
for any monomials $x^{(\be)}, x^{(\ga)}\in \Cal A_q$. Relation (1)
follows easily from relations (i)--(iii) and {\bf 3.1} (i). Relation
(2) is due to the fact that $\Cal K_i$ acts as an algebra
automorphism of $\Cal A_q$.

By {\bf 1.3} (i)--(ii) \& (iii),
we shall prove that the endomorphism $e_i$ $(1\le i<n)$
is a $(\text{id},\sigma_i\sigma_{i+1}^{-1})$-derivation
and $f_i$ $(1\le i<n)$
is a $(\sigma_i^{-1}\sigma_{i+1},\text{id})$-derivation, which implies
relations (3) and (4).

Noting that
$\Th(\vn_{i+1}-\vn_i)=\si_i\,\si_{i+1}$ and
using Proposition 3.1 (6), we see that
$x^{(\vn_i)}\pa_{i+1}$
is a $(\Th(\vn_i-\vn_{i+1})\,\si_{i+1},\si_{i+1}^{-1})$-derivative
of $\Cal A_q$, i.e. $(\si_i^{-1},\si_{i+1}^{-1})$-derivative,
hence, $x^{(\vn_i)}\,\pa_{i+1}\,\si_i$
is a $(\text{id},\si_i\,\si_{i+1}^{-1})$-derivative (since $\si_i$ is an
automorphism of $\Cal A_q$).
On the other hand,
$x^{(\vn_{i+1})}\pa_i$ is a
$(\Th(\vn_{i+1}-\vn_i)\si_i^{-1},\si_i)$-derivative of $\Cal A_q$, i.e.
$(\si_{i+1},\si_i)$-derivative, so
$\si_i^{-1}\,x^{(\vn_{i+1})}\,\pa_i$ is a
$(\si_i^{-1}\,\si_{i+1},\text{id})$-derivative (since $\si_i^{-1}$ is an
automorphism of $\Cal A_q$).

Therefore, we complete the proof of Theorem 4.1.
\hfill\qed
\enddemo

\proclaim{Corollary} \ For any monomial $x^{(\be)}\in \Cal A_q$, set
$$
\gather
e_i(x^{(\be)})=\bigl(\,x_i\pa_{i+1}\sigma_i\,\bigr)\,(x^{(\be)})
                 =[\be_i+1]\,x^{(\be+\vn_i-\vn_{i+1})} \qquad (1\le i<n),\\
f_i(x^{(\be)})=\bigl(\,\sigma_i^{-1}x_{i+1}\pa_i\,\bigr)\,(x^{(\be)})
                 =[\be_{i+1}+1]\,x^{(\be-\vn_i+\vn_{i+1})} \qquad (1\le i<n),\\
k_i(x^{(\be)})=\sigma_i(x^{(\be)}) \qquad (1\le i\le n),\\
k_i^{-1}(x^{(\be)})=\sigma_i^{-1}(x^{(\be)}) \qquad (1\le i\le n).
\endgather
$$
Formulas $(\text{\rm i})$--$(\text{\rm iv})$ define the structure of a $U_q(\frak {gl}_n)$-module
algebra on $\Cal A_q$.
\qed
\endproclaim

\noindent
{\it Remark.} \ According to the interpretation
given in section {\bf 2.4},
the realization above is also valid for the quantum
$n$-space $k[A_q^{n|0}]$,
which can be viewed as an improvement of Hayashi's one
(cf. \cite{5} \& \cite{6}).
It should be noticed that Hayashi's
realization of $U_q(\frak {sl}_n)$ was
carried out over the usual polynomial algebra in $n$ variables
and we can verify that his realization cannot
make the polynomial algebra
a $U_q(\frak {sl}_n)$-module
algebra in the sense of {\bf 1.3} (cf.
\cite{1} \& \cite{16}).

\medskip
\noindent
{\bf 4.2} \ As a direct consequence of
Theorem 4.1 or Corollary 4.1,
we consider the submodule structures of $\Cal A_q$,
where $\Cal A_q=\Cal A_q(n)$
if $\ch(q)=0$, $\Cal A_q=\Cal A_q(n,\bold 1)$ if
$\ch(q)=l\,(\ge 3)$.
Denote $|\,\al\,|:=\sum_{i=1}^n \al_i$ by the degree of
$x^{(\al)}\in\Cal A_q$, set $N:=|\,\tau\,|=n\,(l-1)$. Let
$\Cal A_q^{(s)}:
=\langle \,x^{(\al)}\mid |\,\al\,|=s\,\rangle$.
Then $\Cal A_q(n)=\bigoplus_{s\ge 0}\Cal A_q(n)^{(s)}$, and
 $\Cal A_q(n,\bold 1)=\bigoplus_{s=0}^N
\Cal A_q(n,\bold 1)^{(s)}$ when $\ch(q)=l$ are
$\Bbb Z_+$-graded algebras.

\proclaim{Proposition} \ The subspace $\Cal A_q^{(s)}$ of
homogeneous elements of degree $s$ is a $U_q(\frak {sl}_n)$-submodule of
the quantum divided power (restricted) algebra $\Cal A_q$.

$(1)$ If $\ch(q)=0$, $\Cal A_q(n)^{(s)}$
is generated by the highest weight vector
$x^{(s\,\vn_1)}$ $($where $s\,\vn_1=(s,0,\cdots,0)\,)$, which is
isomorphic to the simple module $V(s\,\lambda_1)$
$($where $\lambda_1=\vn_1$
is the $1$-st fundamental weight of $\frak g\,)$.

$(2)$ If $\ch(q)=l\ge 3$, $\Cal A_q(n,\bold 1)^{(s)}$ is generated by
the highest weight vector
$x^{(\,(l-1)\vn_1+\cdots+(l-1)\vn_{i-1}+s_i\vn_i)}$ $($where
$s=(i-1)(l-1)+s_i$, $0\le s_i\le l-1$ for $1\le i\le n\,)$, which is isomorphic to the simple
module $V(\la)$ $($where $\la=(l-1-s_i)\la_{i-1}+s_i\,\la_i$, $\la_i=\vn_1+
\cdots+\vn_i$
$(1\le i<n)$ is
the $i$-th fundamental weight of $\frak g$ with $\la_0=\la_n=0\,)$.
\endproclaim
\demo{Proof} \
Obviously, we see from formulae {\bf 4.1} (iii)--(v) that the action of
$U_q(\frak {sl}_n)$ on $\Cal A_q$ stabilizes $\Cal A_q^{(s)}$. Moreover,
in the case when $\ch(q)=0$, we note that for $1\le i<n$,
$$
\split
e_i(x^{(s\,\vn_1)})&=0, \\
\Cal K_i(x^{(s\,\vn_1)})&
=q^{\delta_{i,1}s}\,x^{(s\,\vn_1)}\\
&=q^{\langle s\,\la_1,\vn_i-\vn_{i+1}\rangle}
\,x^{(s\,\vn_1)},
\endsplit
$$
which imply the vector $x^{(s\,\vn_1)}$ is a highest weight vector with
highest weight $s\,\la_1$.
Again for any $x^{(\be)}\in\Cal A_q(n)^{(s)}$, $|\be|=s$, set
$s_i=s-\sum_{j\le i}\be_j$ ($1\le i<n$), that is, $\be_1=s-s_1$,
$\be_2=s_1-s_2, \cdots, \be_{n-1}=s_{n-2}-s_{n-1}$, $\be_n=s_{n-1}$,
from formula {\bf 4.1} (v) and by induction on $s_i$, we get
$$
\split
f_{n-1}^{s_{n-1}}&\cdots f_1^{s_1}(x^{(s\,\vn_1)})\\
&=[s_1]!\cdots[s_{n-1}]!\,
x^{((s-s_1)\vn_1+(s_1-s_2)\vn_2+\cdots+(s_{n-2}-s_{n-1})\vn_{n-1}+s_{n-1}\,\vn_n)}\\
&=[s_1]!\cdots[s_{n-1}]!\,x^{(\be)},
\endsplit
$$
where $[m]!=[m]\cdots [2]\cdot [1]$, for $m\in\Bbb Z_+$. Thus we show
that $\Cal A_q(n)^{(s)}\cong V(s\,\la_1)$ is indeed a simple highest
weight module with $s\,\la_1$ as its highest weight, where  $x^{(s\,\vn_1)}$
is its a highest weight vector.

In the case when $\ch(q)=l\,(\ge 3)$, for $0\le s\le N$, there exist $i$
and $s_i$ such that $1\le i\le n$, $s=(i-1)(l-1)+s_i$ with $0\le s_i\le l-1$.
Consider the vector $x^{(\bold s)}\in\Cal A_q(n,\bold 1)^{(s)}$
where $\bold s=(l-1)\vn_1+\cdots+(l-1)\vn_{i-1}+
s_i\vn_i$. From formulae {\bf 4.1} (iii)--(v) and noting
$[\,l\,]=0$, we conclude that for $1\le j<n$,
$$
\split
\Cal K_j(x^{(\bold s)})
&=q^{\langle(l-1-s_i)\la_{i-1}+s_i\la_i, \,\vn_j-\vn_{j+1}\rangle}\,
x^{(\bold s)},\\
e_j(x^{(\bold s)})&=0,
\endsplit
$$
which imply $x^{(\bold s)}$
is a highest weight vector with highest weight
$(l-1-s_i)\la_{i-1}+s_i\la_i$.

On the other hand, for any $x^{(\be)}\in \Cal A_q(n,\bold 1)^{(s)}$,
we have $|\be\,|=s$, $0\le \be_i\le l-1$ ($1\le i\le n$). We denote
$r$ by the last ordinal dumber with $\be_r\ne 0$ for
$n$-tuple $\be=(\be_1,\cdots,\be_n)$. Then $r\ge i$ if $s_i\ne 0$, and
$r\ge i-1$ if $s_i=0$.
Hence, in terms of formula {\bf 4.1} (v), we obtain that

Case (i): if $s_i\ge \be_r\,(>0)$, then
$$
\split
f_{r-1}^{\be_r}\cdots f_i^{\be_r}(x^{(\bold s)})
&=(\,[\be_r]!\,)^{r-i}\,x^{(\bold s-\be_r\vn_i+\be_r\vn_r)}\\
&=(\,[\be_r]!\,)^{r-i}\,x^{(\bold s')}\,x^{(\be_r\vn_r)}\ne 0,
\endsplit
$$
where
$\bold s':=\bold s-\be_r\vn_i=(l-1)\vn_1+\cdots +(l-1)\vn_{i-1}$.

Case (ii): if $\be_r>s_i\,(\ge 0)$, then
$$
\split
f_{r-1}^{\be_r-s_i}&\,f_{r-2}^{\be_r-s_i}\cdots f_{i-1}^{\be_r-s_i}
\,f_{r-1}^{s_i}\cdots f_i^{s_i}\,(x^{(\bold s)})\\
&=(\,[s_i]!\,[\be_r-s_i]!\,)^{r-i}\,[s_i+1][s_i+2]\cdots [\be_r]\,
x^{(\bold s-(\be_r-s_i)\vn_{i-1}-s_i\vn_i+\be_r\vn_r)}\\
&=(\,[s_i]!\,[\be_r-s_i]!\,)^{r-i}\,[s_i+1][s_i+2]\cdots [\be_r]\,
x^{(\bold s')}\,x^{(\be_r\vn_r)}\ne 0,
\endsplit
$$
where
$\bold s':=\bold s-(\be_r-\vn_i)\vn_{i-1}-s_i\vn_i=(l-1)\vn_1+\cdots
((l-1)\vn_{i-2}+(l-1-\be_r+s_i)\vn_{i-1}$.

Set $\be':=\be-\be_r\vn_r$ and use an induction on $\bold s$. At first,
the argument holds for $\bold s=\vn_1=\la_1$ (see the following
Remark). Assume that there exists a word $\omega$ in $U_q(\frak {sl}_n)$
constructed by some suitable $f_j$'s (where $j<r$) such that
$\omega\,(x^{(\bold s')})=c\,x^{(\be')}$ ($c\in k^*$).
Note that $f_j\,(x^{(\be_r\vn_r)})=0$ for those $f_j$ $(j<r)$. Thus
we get
$$
\split
\omega\,(x^{(\bold s')}\,x^{(\be_r\vn_r)})
&=\omega\,(x^{(\bold s')})\,x^{(\be_r\vn_r)}\\
&=c\,x^{(\be')}\,x^{(\be_r\vn_r)}=c\,x^{(\be)}\ne 0.
\endsplit
$$
Since $x^{(\be)}$ (with $|\be|=s$) is arbitrary,
$\Cal A_q(n,\bold 1)^{(\bold s)}$ is an indecomposable module
generated by the highest weight vector $x^{(\bold s)}$.

Finally, by virtue of {\bf 4.1} (v), we find both $e_i$'s and
$f_i$'s act nilpotently on $\Cal A_q(n,\bold 1)$, namely,
$e_i^l|_{\Cal A_q(n,\bold 1)}\equiv 0 \equiv f_i^l|_{\Cal
A_q(n,\bold 1)}$. This implies the $U_q(\frak {sl}_n)$-module $\Cal
A_q(n,\bold 1)$ is completely reductive. Consequently, we derive
that $\Cal A_q(n,\bold 1)^{(s)}\cong V((l-1-s_i)\la_{i-1}+s_i\la_i)$
is a simple highest weight module. \hfill\qed
\enddemo

\noindent
{\it Remark.} \ Particularly, $\langle x_1,\cdots,x_n\rangle\cong V(\la_1)$
with $e_i(x_j)=\de_{i+1,j}\,x_i$, $f_i(x_j)=\de_{ij}\,x_{i+1}$ and
$\Cal K_i(x_1)=q^{\langle \la_1, \vn_i-\vn_{i+1}\rangle}\,x_1$.
Observe that
the conclusion (1) of Proposition 4.2 is valid for the quantum $n$-space
$k[A_q^{n|0}]$ when $\ch(q)=0$ (see Remark 4.1).

\medskip
\noindent
{\bf 4.3} \ As the dual object of the quantum $n$-space (cf. \cite{11}),
we can consider its submodule structures of
the quantum exterior algebra
$k[A_q^{0|n}]=\Lambda_q(n)=\bigoplus_{s=0}^n\La_q(n)_{(s)}$ where
$\Lambda_q(n)_{(s)}:=
\langle\,x_{i_1}\cdots x_{i_s}\mid 1\le i_1<\cdots<i_s\le n\,\rangle$.
The known fact below is independent of $\ch(q)$.
\proclaim{Proposition} \ The subspace $\Lambda_q(n)_{(s)}$ of
homogeneous elements of degree $s$ is a $U_q(\frak {sl}_n)$-submodule of
$\Lambda_q(n)$. It is generated by the highest weight
vector $x_1\cdots x_s$ and is isomorphic to the simple module
$V(\lambda_s)$ $($where $\lambda_s=\vn_1+\cdots+\vn_s$ is the $s$-th
fundamental weight of $\frak {sl}_n)$.
In other words, the quantum exterior
algebra $\Lambda_q(n)$ is the direct sum of
those all basic simple modules of
$U_q(\frak {sl}_n)$.
\endproclaim
\demo{Proof} \ We can identify elements in $\La_q(n)_{(s)}$ with those
in tensor algebra $\Cal T(V)$ with a $U_q(\frak {sl}_n)$-action induced from the
action on $V$ via the $s$-th
comultiplication $\De^{(s)}$, where $V=\langle x_1,\cdots,x_n\rangle$ is
the first basic module (see Remark 4.2). Noting $x_i^2=0$ in $\La_q(n)$,
we readily obtain that for $1\le i<n$,
$$
\split
\Cal K_i(x_1\cdots x_s)&=q^{\delta_{is}}\;x_1\cdots x_s\\
        &=q^{\langle\lambda_s,\vn_i-\vn_{i+1}\rangle}\;x_1\cdots x_s,\\
e_i(x_1\cdots x_s)&=0 \qquad (\text{\it since \;} \;x_i^2=0),\\
f_{j_s-1}\cdots f_s(x_1\cdots x_s)&=x_1\cdots x_{s-1}x_{j_s} \qquad
(j_s\ge s),
\endsplit
$$
thus for any $1\le j_1<\cdots<j_s\le n$, one has
$$
(f_{{j_1}-1}\cdots f_1)\cdots (f_{{j_s}-1}\cdots f_s)\;(x_1\cdots x_s)
=x_{j_1}\cdots x_{j_s}.
$$
Consequently, $x_1\cdots x_s$ is a highest weight vector of weight
$\lambda_s$ and generates the simple submodule $\Lambda_q(n)_{(s)}$.
\hfill\qed
\enddemo

\noindent
{\bf 4.4} \ Theorem 4.1 has indicated the presentation in $\Cal W_q(2n)$
of the generators of $U_q(\frak {sl}_n)$. Actually, we can describe
explicitly the presentation of all ``{\it roots vectors}" of $U_q(\frak {sl}_n)$
under our realization, which is coincident with one of four kinds of
roots vectors introduced by G. Lusztig in terms of its braid
(automorphism) group action (cf. Lemma 39.3.2, Corollary 40.2.2 in
\cite{9}). To do so, we need some notions.

Recall some known facts about $U_q(\frak {sl}_n)$. Let $\al_i=\vn_i-\vn_{i+1}$
$(1\le i<n$) be the simple roots of $\frak {sl}_n$ (see section {\bf 1.4}),
$s_i$ $(1\le i<n)$ the reflections determined by those $\al_i$ respectively and
generate the Weyl group $W=S_n$ (the permutation group). The general
positive roots in $\De^+$ take the form $\al_{ij}=\al_i+\al_{i+1}+\cdots+\al_{j-1}
=\vn_i-\vn_j$ ($1\le i<j\le n$) with $\al_{i,i+1}=\al_i$. Fix a reduced
element $\omega_0=s_{i_1}s_{i_2}\cdots s_{i_N}$ of maximal
length in $W$. Then $\al=s_{i_1}s_{i_2}\cdots s_{i_{p-1}}(\al_{i_p})$ as
$p=1,2,\cdots,N$ runs through all the positive roots of $\De^+$

Let $T_i$ be the braid automorphism of $U_q(\frak  {sl}_n)$ corresponding to
$s_i$ ($1\le i<n$), which have been introduced as $T_{i,-1}''$ by Lusztig
in \S37.1.3 \cite{9} and take the form
$$
\gather
T_i(\Cal K_\mu)=\Cal K_{s_i(\mu)}, \qquad
T_i(e_i)=-f_i\,\Cal K_i^{-1}, \qquad T_i(f_i)=-\Cal K_i\,e_i;\tag{1}\\
T_i(e_j)=e_j, \qquad T_i(f_j)=f_j, \qquad (|i-j|>1);\\
T_i(e_j)=e_i\,e_j-q\,e_j\,e_i, \qquad (|i-j|=1);\\
T_i(f_j)=f_jf_i-q^{-1}\,f_i\,f_j, \qquad (|i-j|=1)\\
\endgather
$$
where $e_i$ ($1\le i<n$) are the simple roots vectors of $U_q(\frak {sl}_n)$
corresponding to
the simple roots $\al_i$ respectively, and $f_i$ ($1\le i<n$) the negative
simple roots vectors to those $-\al_i$. As we know, the elements below
$$
e_{\al}=T_{i_1}\cdots T_{i_{p-1}}(e_{i_p}), \qquad
f_{\al}=T_{i_1}\cdots T_{i_{p-1}}(f_{i_p})
$$
are called roots vectors of $U_q(\frak {sl}_n)$ associated to boots $\pm\al
=\pm s_{i_1}s_{i_2}\cdots s_{i_{p-1}}(\al_{i_p})$ ($p=1, 2, \cdots, N$)
respectively, where $\omega_0=s_{i_1}s_{i_2}\cdots s_{i_N}$
is as above. Normally, for the Weyl group of $\frak {sl}_n$, we can take such
a reduced element $\omega_0$ as
$$
\omega_0=s_1s_2s_1s_3s_2s_1\cdots s_{n-2}s_{n-3}\cdots s_2s_1s_{n-1}s_{n-2}\cdots s_2s_1,
\leqno(2)
$$
which gives rise to a well-known ordering of $\De^+$ as
$$
\al_{12}, \;\al_{13}, \;\al_{23}, \;\al_{14}, \;\al_{24}, \;
\al_{34}, \;\cdots, \;
\al_{1n}, \;\al_{2n}, \;\cdots, \;\al_{n-1,n}.\leqno(3)
$$
Hence, all the positive roots vectors of $U_q(\frak {sl}_n)$ associated
to the ordering (3) of $\De^+$ are as follows.
$$
\gather
e_{\al_{12}}=e_1, \tag{4}\\
e_{\al_{13}}=T_1(e_2), \qquad
e_{\al_{23}}=T_1T_2(e_1)=e_2,\\
e_{\al_{14}}=T_1T_2T_1(e_3), \qquad
e_{\al_{24}}=T_1T_2T_1T_3(e_2), \qquad
e_{\al_{34}}=T_1T_2T_1T_3T_2(e_1)=e_3, \\
\cdots\cdots\\
e_{\al_{1n}}=T_1T_2T_1T_3T_2T_1\cdots T_{n-2}T_{n-3}\cdots T_2T_1(e_{n-1}),\\
e_{\al_{2n}}=T_1T_2T_1T_3T_2T_1\cdots T_{n-2}T_{n-3}\cdots T_2T_1T_{n-1}(e_{n-2}),\\
\cdots\cdots\\
e_{\al_{n-1,n}}=T_1T_2T_1T_3T_2T_1\cdots T_{n-2}T_{n-3}\cdots T_2T_1T_{n-1}T_{n-2}\cdots
T_2(e_1)=e_{n-1}.
\endgather
$$

\noindent
{\bf 4.5} \
We introduce here some $q$-differential operators in $\Cal W_q(2n)$.
Set $E_{ij}:=x_i\pa_j$, for any $1\le i, \;j\le n$. Denote
$e_{ij}:=E_{ij}\,\si_i$ for $1\le i<j\le n$ and
$e_{ij}:=\si_j^{-1}\,E_{ij}$ for $1\le j<i\le n$.
Since
$$
\split
E_{ij}(x^{(\be)})
&=\bigl(\,x_i\pa_j\,\bigr)\,(x^{(\be)})\\
&=[\be_i+1]\,q^{\vn_i*(\be-\vn_j)-\vn_j*\be}\,x^{(\be+\vn_i-\vn_j)}\\
&=[\be_i+1]\,q^{(\vn_i-\vn_j)*\be-\vn_i*\vn_j}\,x^{(\be+\vn_i-\vn_j)}\\
&=[\be_i+1]\,q^{-\sum_{i\le s<j}\be_s}\,x^{(\be+\vn_i-\vn_j)},\qquad (i<j)\\
E_{ij}(x^{(\be)})
&=\bigl(\,x_i\pa_j\,\bigr)\,(x^{(\be)})\\
&=[\be_i+1]\,q^{(\vn_i-\vn_j)*\be-\vn_i*\vn_j}\,x^{(\be-\vn_j+\vn_i)}\\
&=[\be_i+1]\,q^{\sum_{j\le s<i}\be_s-1}\,x^{(\be-\vn_j+\vn_i)},\qquad (i>j)
\endsplit
$$
thus we get
$$\gather
e_{ij}(x^{(\be)})
=(E_{ij}\sigma_i)\,(x^{(\be)})
=[\be_i+1]\,q^{-\sum_{i<s<j}\be_s}\,x^{(\be+\vn_i-\vn_j)},\quad (i<j)\tag{1}\\
e_{ij}(x^{(\be)})
=(\sigma_j^{-1}E_{ij})\,(x^{(\be)})
=[\be_i+1]\,q^{\sum_{i>s>j}\be_s}\,x^{(\be-\vn_j+\vn_i)}.\quad (i>j)\tag{2}
\endgather
$$

\proclaim{Lemma} \ $(\text{\rm i})$ If $i<j$, then for any $i<k<j$, we have
$e_{ij}=e_{ik}\,e_{kj}-q\,e_{kj}\,e_{ik}$;

$(\text{\rm ii})$ If $i>j$, then for any $i>k>j$, we have
$e_{ij}=e_{ik}\,e_{kj}-q^{-1}\,e_{kj}\,e_{ik}$;

$(\text{\rm iii})$ For $i<j$, we have $e_{ij}\,e_{ji}-e_{ji}\,e_{ij}
=\frac{\si_i\,\si_j^{-1}-\si_i^{-1}\,\si_j}{q-q^{-1}}$,
where $\si_i\,\si_j^{-1}=\Cal K_{\vn_i-\vn_j}$.
\endproclaim
\demo{Proof}
(i) For $i<k<j$. Since
$[m+1]-q\,[m]=q^{-m}$, it follows from (1) that
$$
\split
&\bigl(\,e_{ik}\,e_{kj}-q\,e_{kj}\,e_{ik}\,\bigr)
  \,(x^{(\be)})\\
&\qquad=[\be_k+1]\,q^{-\sum_{k<s<j}\be_s}\,e_{ik}(x^{(\be+\vn_k-\vn_j)})
  -q\,[\be_i+1]\,q^{-\sum_{i<s<k}\be_s}\,e_{kj}(x^{(\be+\vn_i-\vn_k)})\\
&\qquad=(\,[\be_k+1]\,[\be_i+1]-q\,[\be_i+1]\,[\be_k]\,)\,q^{-\sum_{i<s<j}\be_s+\be_k}\,
  x^{(\be+\vn_i-\vn_j)}\\
&\qquad=[\be_i+1]\,q^{-\sum_{i<s<j}\be_s}\,x^{(\be+\vn_i-\vn_j)}
=e_{ij}\,(x^{(\be)}), \qquad \forall \; x^{(\be)}\in\Cal A_q.
\endsplit
$$

(ii) For $i>k>j$.
Since
$[m+1]-q^{-1}\,[m]=q^m$, it follows from (2) that
$$
\split
&\bigl(\,e_{ik}\,e_{kj}-q^{-1}\,e_{kj}\,e_{ik}\,\bigr)
  \,(x^{(\be)})\\
&\qquad=[\be_k+1]\,q^{\sum_{k>s>j}\be_s}\,e_{ik}(x^{(\be+\vn_k-\vn_j)})
  -q^{-1}\,[\be_i+1]\,q^{\sum_{i>s>k}\be_s}\,e_{kj}(x^{(\be+\vn_i-\vn_k)})\\
&\qquad=(\,[\be_k+1]\,[\be_i+1]-q^{-1}\,[\be_i+1]\,[\be_k]\,)\,q^{\sum_{i>s>j}\be_s-\be_k}\,
  x^{(\be+\vn_i-\vn_j)}\\
&\qquad=[\be_i+1]\,q^{\sum_{i>s>j}\be_s}\,x^{(\be+\vn_i-\vn_j)}
=e_{ij}\,(x^{(\be)}), \qquad \forall \; x^{(\be)}\in\Cal A_q.
\endsplit
$$

(iii) For $i<j$, by (1) \& (2), we get
$$
\split
\bigl(\,e_{ij}&\, e_{ji}-e_{ji}\,e_{ij}\,\bigr)\,
(x^{(\be)})\\
&=(\,[\be_i]\,[\be_j+1]-[\be_j]\,[\be_i+1]\,)\,x^{(\be)}\\
&=[\be_i-\be_j]\,x^{(\be)}\\
&=\frac{\si_i\,\si_j^{-1}-\si_i^{-1}\,\si_j}{q-q^{-1}}\,(x^{(\be)}),
\qquad \forall \; x^{(\be)}\in \Cal A_q.
\endsplit
$$

Thus we conclude the proof.
\qed
\enddemo

We introduce two $q$-brackets as follows.
$$\gather
[e_{ik},e_{kj}]_q:=
e_{ik}\,e_{kj}-q\,e_{kj}\,e_{ik},\qquad\qquad\text{\it for \;} \;i<k<j.\tag{3}\\
[e_{ik},e_{kj}]_{q^{-1}}:=
e_{ik}\,e_{kj}-q^{-1}\,e_{kj}\,e_{ik}. \qquad \text{\it for \;} \; i>k>j.
\tag{4}
\endgather
$$
So we arrive at $e_{ij}=[e_{ik},e_{kj}]_q$ for $i<j$, and
$e_{ij}=[e_{ik},e_{kj}]_{q^{-1}}$ for $i>j$.

On the other hand, the preceding Lemma also indicates
$e_{ij}$ expressed by the $q$-brackets in two cases are independent of
the choice of $k$.
Hence, we deduce that
$$
\split
e_{ij}
&=[\cdots[e_{i,i+1},e_{i+1,i+2}]_q,\cdots,e_{j-2,j-1}]_q,e_{j-1,j}]_q\\
&=[e_{i,i+1},[e_{i+1,i+2},\cdots,[e_{j-2,j-1},e_{j-1,j}]_q\cdots]_q,\qquad\qquad (i<j)\\
e_{ij}
&=[\cdots[e_{i,i-1},e_{i-1,i-2}]_{q^{-1}},\cdots,e_{j+2,j+1}]_{q^{-1}},e_{j+1,j}]_{q^{-1}}\\
&=[e_{i,i-1},[e_{i-1,i-2},\cdots,[e_{j+2,j+1},e_{j+1,j}]_{q^{-1}}\cdots]_{q^{-1}},\qquad (i>j)
\endsplit\tag{5}
$$

\noindent
{\it Remark.} \
It should be mentioned that the validity of $q$-bracket $[\,,\,]_q$ defined in (3) (resp.
$[\,,\,]_{q^{-1}}$ in (4)\,) only involves one direction of an ordering of
$e_{ij}$, which is a different point from the defining relation of
$T_i$ on $e_j$ in {\bf 4.4} (1), however, it doesn't affect our conclusion
of the next Proposition 4.6.

\medskip
\noindent
{\bf 4.6} \
We now shall show that the positive roots vectors $e_{\al_{ij}}$ of
$U_q(\frak {sl}_n)$ given in {\bf 4.4} (4) just correspond to
those $q$-differential operators $e_{ij}$ described in {\bf 4.5} (1)
in the sense of our
realization (cf. Theorem 4.1), that is, if we identify $e_i$
(resp. $f_i$) with $e_{i,i+1}$ (resp. $e_{i+1,i}$), as well as
$\Cal K_i$ with $\si_i\,\si_{i+1}^{-1}$ for $1\le i<n$,
then all positive roots vectors $e_{\al_{ij}}$ can be identified with
the above $q$-differential operators $e_{ij}$. More generally, we have

\proclaim{Proposition} \ Identifying the generators of $U_q(\frak {sl}_n)$
with the certain $q$-differential operators in $\Cal W_q(2n)$, i.e.
$e_i:=e_{i,i+1}$, $f_i:=e_{i+1,i}$, $\Cal K_i:=\si_i\,\si_{i+1}^{-1}$
with $1\le i<n$. Then we have

$(\text{\rm i})$ \ $e_{ij}$ $(i<j)$ correspond to the positive root vectors
$e_{\al_{ij}}$ associated to those positive roots $\al_{ij}=\vn_i-\vn_j$
$(i<j)$, i.e. $e_{\al_{ij}}:=e_{ij}$.

$(\text{\rm ii})$ \ $e_{ij}$ $(i>j)$ correspond to the negative root vectors
$f_{\al_{ji}}$ associated to the negative roots $-\al_{ji}=\vn_i-\vn_j$
$(i>j)$, i.e. $f_{\al_{ji}}:=e_{ij}$. \hfill\qed
\endproclaim

Obviously, it suffices to prove the first claim. For this purpose,
according to the ordering of $\De^+$ made in {\bf 4.4} (3),
we will use an induction on the length of ordered subword of $\omega_0$ as
in {\bf 4.4} (2)
to show the required identification relation of all positive roots
vectors $e_{\al_{ij}}$ ($i<j$).

To do so, we will first establish an
auxiliary Lemma and make use of the
following facts due to Lusztig (cf. \S39.2.4 in \cite{9}) in our
argument.
$$
\gather
T_iT_jT_i=T_jT_iT_j, \qquad |i-j|=1,\tag{1}\\
T_iT_j(e_i)=e_j, \qquad |i-j|=1,\tag{2}\\
T_i(e_j)=e_j, \qquad |i-j|>1.\tag{3}
\endgather
$$

\proclaim{Lemma} \ With the identification as the preceding
Proposition, we have

$(\text{\rm i})$ \ If $i+1<j$, then
$[e_{ij}, T_i(e_i)]_q=q\,f_i\,\Cal K_i^{-1}\,e_{ij}
-e_{ij}\,f_i\,\Cal K_i^{-1}=e_{i+1,j}$.

$(\text{\rm ii})$ If $i+1<j$, then
$[T_i(f_i), e_{ji}]_{q^{-1}}=q^{-1}\,e_{ji}\,\Cal K_ie_i-\Cal K_ie_i\,e_{ji}
=e_{j,i+1}$.
\endproclaim
\demo{Proof}
(i) When $i+1<j$, $\forall \; x^{(\be)}\in\Cal A_q$, using {\bf 4.4} (1),
{\bf 4.5} (1) \& (2), we have
$$
\split
\bigl(\,q\,f_i\,\Cal K_i^{-1}\,e_{ij}&
-e_{ij}\,f_i\,\Cal K_i^{-1}\,\bigr)
(x^{(\be)})\\
&=q^{1-\sum_{i<s<j}\be_s}[\be_i+1]\,f_i\,\Cal K_i^{-1}
(x^{(\be+\vn_i-\vn_j)})\\
&\quad \, -q^{-\be_i+\be_{i+1}}\,[\be_{i+1}+1]\,e_{ij}
(x^{(\be-\vn_i+\vn_{i+1})})\\
&=q^{-\sum_{i\le s<j}\be_s+\be_{i+1}}\,[\be_i+1]\,[\be_{i+1}+1]\,
x^{(\be+\vn_{i+1}-\vn_j)}\\
&\quad \, -q^{-\be_i+\be_{i+1}-\sum_{i<s<j}\be_s-1}\,
[\be_{i+1}+1]\,[\be_i]\,
x^{(\be+\vn_{i+1}-\vn_j)}\\
&=q^{-\be_i-\sum_{i+1<s<j}\be_s}\,[\be_{i+1}+1]\,(\,[\be_i+1]
-q^{-1}\,[\be_i]\,)\,x^{(\be+\vn_{i+1}-\vn_j)}\\
&=q^{-\sum_{i+1<s<j}\be_s}\,
[\be_{i+1}+1]\,x^{(\be+\vn_{i+1}-\vn_j)}
=e_{i+1,j}(x^{(\be)}),
\endsplit
$$
so the conclusion is true.

(ii) When $j>i+1$, $\forall \; x^{(\be)}\in\Cal A_q$, using {\bf 4.4} (1),
{\bf 4.5} (1) \& (2), we have
$$
\split
\bigl(\,q^{-1}\,e_{ji}\,\Cal K_i\,e_i&-\Cal K_i\,e_i\,e_{ji}\,\bigr)
(x^{(\be)})\\
&=\bigl(\,q^{-1+\be_i-\be_{i+1}+1+\sum_{j>s>i}\be_s}\,[\be_i+1]\,[\be_j+1]\bigr.\\
&\quad \,\bigl.-q^{\sum_{j>s>i}\be_s+\be_i-\be_{i+1}+1}\,[\be_i]\,[\be_j+1]\,\bigr)
x^{(\be+\vn_j-\vn_{i+1})}\\
&=q^{\sum_{j>s>i+1}\be_s}\,[\be_j+1]\,x^{(\be+\vn_j-\vn_{i+1})}
=e_{j,i+1}(x^{(\be)}),
\endsplit
$$
so the claim is true.
\hfill\qed
\enddemo

We are now in the position to show Proposition 4.6.

\demo{Proof of Proposition 4.6} \ We suffice to verify our identification
$e_{\al_{ij}}:=e_{ij}$ for the case when $i<j$. We adopt the ordering of
positive roots vectors $e_{\al_{ij}}$ as in the list {\bf 4.4} (4) and use
an induction on $k$ where $1\le k< n$.

For $k=1$, this is just
the identification $e_{\al_1}=e_1:=e_{12}$.

For $k=2$,
by the identification and using {\bf 4.5} (3) \& {\bf 4.6} (2),
we have
$$\split
e_{\al_{13}}&=T_1(e_2)=e_1\,e_2-q\,e_2\,e_1:=
e_{12}\,e_{23}-q\,e_{23}\,e_{12}\\
&=[e_{12},e_{23}]_q=e_{13},\\
e_{\al_{23}}&=T_1T_2(e_1)=e_2:=e_{23}.
\endsplit
$$

Assume that the claim is valid
for the case $U_q(\frak {sl}_{n-1})$.
Thus for the case $U_q(\frak {sl}_n)$,
we further need to
show that in the list {\bf 4.4} (4) there hold
$$
e_{\al_{1n}}:=e_{1n},
\quad e_{\al_{2n}}:=e_{2n}, \; \cdots, \quad
e_{\al_{n-1,n}}:=e_{n-1,n}.
$$

At first, according to the above assumption, we have
$$
e_{\al_{1,n-1}}=
T_1T_2T_1T_3T_2T_1\cdots
T_{n-4}\cdots T_2T_1T_{n-3}\cdots T_2T_1
(e_{n-2}):=e_{1,n-1}.
$$
Thus, noting (3) and using {\bf 4.5} (3), we have
$$
\split
e_{\al_{1n}}&=
T_1T_2T_1T_3T_2T_1\cdots T_2T_1T_{n-3}\cdots T_2T_1T_{n-2}
\cdots T_2T_1(e_{n-1})\\
&=T_1T_2T_1T_3T_2T_1\cdots T_2T_1T_{n-3}\cdots T_2T_1T_{n-2}
(e_{n-1})\\
&=T_1T_2T_1T_3T_2T_1\cdots T_{n-3}\cdots T_2T_1\,(e_{n-2}\,
e_{n-1}-
q\,e_{n-1}\,e_{n-2}\,)\\
&=T_1T_2T_1T_3T_2T_1\cdots T_{n-3}\cdots T_2T_1\,(e_{n-2})\cdot
e_{n-1}\\
&\qquad - q\,e_{n-1}\cdot T_1T_2T_1T_3T_2T_1\cdots T_2T_1T_{n-3}
\cdots T_2T_1\,(e_{n-2})\\
&:=e_{1,n-1}\,e_{n-1,n}-q\,e_{n-1,n}\,e_{1,n-1}\\
&=[e_{1,n-1},\,e_{n-1,n}]_q\\
&=e_{1n}.
\endsplit
$$

Let us further suppose that we have proved the identification
$$
e_{\al_{jn}}=T_1T_2T_1T_3T_2T_1\cdots T_{n-2}T_{n-3}\cdots
T_2T_1T_{n-1}T_{n-2}\cdots T_{n-j+1}(e_{n-j})
:=e_{jn},
$$
for $1\le j<n-1$.
Now we want to show $e_{\al_{j+1,n}}:=e_{j+1,n}$.
Note that
$$
\split
T_1T_2T_1T_3T_2T_1\cdots T_j\cdots T_2(f_1)&=f_j,\\
T_1T_2T_1T_3T_2T_1\cdots T_j\cdots T_2(\Cal K_1)&=\Cal K_j,\\
T_1T_2T_1T_3T_2T_1\cdots T_j\cdots T_2T_1(e_1)
&=T_1T_2T_1T_3T_2T_1\cdots T_j\cdots T_2(-f_1\,\Cal K_1^{-1})\\
&=-f_j\,\Cal K_j^{-1}=T_j(e_j).
\endsplit
$$
Using {\bf 4.6} (2) and Lemma 4.6 (i),
we get
$$
\split
e_{\al_{j+1,n}}&=
T_1T_2T_1T_3T_2T_1\cdots T_{n-2}T_{n-3}\cdots
T_2T_1T_{n-1}T_{n-2}\cdots T_{n-j+1}T_{n-j}(e_{n-j-1})\\
&=T_1T_2T_1T_3T_2T_1\cdots T_{n-2}T_{n-3}\cdots
T_2T_1T_{n-1}T_{n-2}\cdots T_{n-j+1}(\,[e_{n-j},e_{n-j-1}]_q)\\
&=[\,T_1T_2T_1T_3T_2T_1\cdots T_{n-2}T_{n-3}\cdots
T_2T_1T_{n-1}T_{n-2}\cdots T_{n-j+1}(e_{n-j}), \\
&\qquad T_1T_2T_1T_3T_2T_1\cdots T_{n-2}T_{n-3}\cdots
T_2T_1T_{n-1}T_{n-2}\cdots T_{n-j+1}(e_{n-j-1})\,]_q\\
&:=[\,e_{jn}, T_1T_2T_1\cdots T_{n-3}\cdots
T_2T_1T_{n-2}T_{n-3}\cdots T_{n-j}T_{n-j-1}T_{n-j-2}(e_{n-j-1})\,]_q\\
&=[\,e_{jn}, T_1T_2T_1T_3T_2T_1\cdots T_{n-3}T_{n-4}\cdots
T_{n-j-1}T_{n-j-2}T_{n-j-3}(e_{n-j-2})\,]_q\\
&\qquad\qquad\cdots\cdots\\
&=[\,e_{jn}, T_1T_2T_1T_3T_2T_1\cdots T_{j+1}T_j\cdots
T_2T_1(e_2)\,]_q\\
&=[\,e_{jn}, T_1T_2T_1T_3T_2T_1\cdots T_jT_{j-1}\cdots T_2T_1(e_1)\,]_q\\
&=[\,e_{jn}, T_j(e_j)\,]_q\\
&=e_{j+1,n}.
\endsplit
$$

Therefore, we complete the proof.
\hfill\qed
\enddemo

\vskip2.5cm

\heading{Acknowledgements}\endheading
\vskip0.2cm

The author would like to express his gratitude
to the National Natural Science Foundation of China (Grant No. 19731004),
the Science Foundation of the University Doctoral Program CNCE,
the Ministry of Education of China and the Shanghai Scientific and
Technical Commission, as well as the Institute of Mathematics of
ECNU, for their supports.
Finally, the author gratefully acknowledges the referee for his/her
very helpful comments so that a couple of paragraphs' expositions
are made clearer, especially, two gaps originally appeared
in Theorem 2.3 and Proposition 2.4 are removed.

\vskip0.20cm

\noindent
\smc Department of Mathematics, East China Normal University, Shanghai 200062,
China

\vskip0.12cm \noindent {\it E-mail address}: nhhu\@math.ecnu.edu.cn

\bye